\documentclass[11pt]{amsart}
\usepackage{amssymb}
\pdfoutput=1

 \title[Contact Diffeomorphisms]{The Group of Contact Diffeomorphisms for
 Compact Contact Manifolds }

\date{June 30, 2010}

 \author{J. Bland  \and T. Duchamp}
	\address{John Bland\\ University of Toronto}
	\address{Tom Duchamp\\ University of Washington}

\email{bland@math.toronto.edu}
\email{duchamp@math.washington.edu}

\thanks{The first author was partially supported by an NSERC grant. The
  second author was partially supported by an NSF grant.}



\numberwithin{equation}{subsection}

\newtheorem{lemma}[equation]{Lemma}
\newtheorem{proposition}[equation]{Proposition}
\newtheorem{theorem}[equation]{Theorem}
\newtheorem{corollary}[equation]{Corollary}

\theoremstyle{definition}

\newtheorem{definition}[equation]{Definition}

\newtheorem{remark}[equation]{Remark}

\newtheorem{examples}[equation]{Examples}
\newcommand\Id{\text{Id}}

\newcommand\jH[1]{{j^{#1}_{H}}}

 \newcommand\Wk[1]{\Gamma^{#1}}
 
 \newcommand\Wka[1]{\Gamma^{{#1},\alpha}}
 \newcommand\Wkcont[1]{\Gamma_{cont}^{#1}(TM)}

\newcommand\WkD[1]{\Gamma^{#1}_{\mathcal{D}_{cont}}}
\newcommand\piD{\pi_{\mathcal{D}}}

\newcommand\piH{\pi_{H}}
\newcommand\dH{d_{H}}
\newcommand\dQ{d_{R}}
\newcommand\deltaQ{\delta_{R}}
\newcommand\DeltaQ{\Delta_{R}}
\newcommand\piQ{\pi_{R}}
\newcommand\HQ{{H_{R}}}
\newcommand\DQ{D_{R}}
\newcommand\DQs{\DQ^{*}}
\newcommand\GQ{G_{R}}
\newcommand\RQ{R}
\newcommand\cRQ{\mathcal{R}}

\newcommand\diff[1]{\mathcal{D}^{#1}(M)}
\newcommand\cdiff[1]{\mathcal{D}_{cont}^{#1}(M)}
\newcommand\cdiffD[1]{\mathcal{D}_{cont}^{#1}(D,\tilde{D})}

\newcommand\cchi[1]{{\chi^{#1}_{cont}}}

\newcommand\locdiff[1]{\mathcal{D}^{s}}

\newcommand\st{\; : \;}

\newcommand\image{\mbox{\rm range }}

\newcommand\inner[1]{\left(#1\right)}

\newcommand\norm[1]{\| #1\|}

\newcommand\pd[2]{\frac{\partial #1 }{ \partial #2}}
\newcommand\db{\bar{\partial}}
\newcommand\bdb{\db_{b}}

\renewcommand\a{\alpha}
\renewcommand\b{\beta}

\newcommand\g{\gamma}

\newcommand\s{\sigma}

\newcommand\n{\eta}

\newcommand\w{\omega}
\newcommand\W{\Omega}

\newcommand\R{\mathbb{R}}
\newcommand\C{\mathbb{C}}

 \newcommand\Lie{\mathcal{L}}

 \newcommand\cD{\mathcal{D}}

 \newcommand\cH{\mathcal{H}}

\newcommand\cO{\mathcal{O}}

 \newcommand\cU{\mathcal{U}}
 \newcommand\cV{\mathcal{V}}

 \newcommand\inter{\mbox{\rule{2.0mm}{0.1mm}\kern0.0mm\rule{0.1mm}{3.0mm} }\,}

 \newcommand\lang{\ {{\vrule width .15in height .03pt
 depth .03pt}\!\!\!\!\!\joinrel{\raise3pt\hbox{$\backslash$}}}\ }

\newcommand{\dist}{\text{dist}}
\newcommand{\XH}{X_{H}}
\newcommand\const[2]{C}

\newcommand{\FX}{F_{X}}            

\newcommand{\Quad}{{Q}}   

\newcommand{\piperp}{\pi^{\perp}}

\newcommand{\modeta}{\;\text{mod }\n}

\newcommand{\apriori}{{\it a priori\/}}

\newcommand\range{\text{\rm range}}

\begin{document}
 \begin{abstract}

   For a compact contact manifold $M^{2n+1}$, it is shown that the
   anisotropic Folland-Stein function spaces $\Gamma^s(M), s \ge (2n+4)$
   form an algebra. The notion of anisotropic regularity is extended to
   define the space of $\Gamma^s$-contact diffeomorphisms, which is shown
   to be a topological group under composition and a smooth Hilbert
   manifold. These results are used in a subsequent paper to analyse the
   action of the group of contact diffeomorphisms on the space of
   CR~structures on a compact, three dimensional manifold.
\end{abstract}

 \maketitle

 \section{Introduction}
\label{sec:intro}

Contact manifolds arise naturally in complex and CR geometry.  The boundary
of a strongly pseudoconvex domain is a contact manifold, and more
generally, any strongly pseudoconvex CR manifold is a contact manifold. In
each case, the $\bdb$--operator, which may be thought to embody the
tangential Cauchy Riemann equations, is a natural operator that arises in
analysis. The associated second order operator $\Box_b$ is anisotropic,
being second order in the holomorphic tangential directions, and only first
order in the transverse directions. In \cite{FS}, Folland and Stein
introduced some anisotropic function spaces, the anisotropic Sobolev spaces
$\Wk{s}$ and the anisotropic Banach spaces $\Wka {s}$ to reflect this
behaviour, and showed that these operators are solvable with good estimates
in these spaces.

In recent years, much attention has been focused on the space of CR
structures which a given compact manifold admits. A theorem of Gray
\cite{gray:1959a} states that all contact structures in the same homotopy
class are equivalent. It is natural, therefore, to fix an underlying
contact structure and study the action of the space of contact
diffeomorphisms on the space of CR~structures that are compatible with the
fixed contact structure.

In \cite{CL}, Cheng and Lee constructed a transverse slice for the action of
the contact diffeomorphism group.  They avoided using the anisotropic
spaces in \cite{CL} by working in the Nash Moser category.  In \cite{B}, we
restricted our attention to the case of the standard $S^3 \subset \C^2$,
and used explicit information to construct an anisotropic Hilbert space
structure on contact diffeomorphisms near the identity.  In this paper, we
show how to generalize the construction in \cite{B}, and show that, for
general compact contact manifolds, contact diffeomorphisms form a smooth
Hilbert manifold modelled on the anisotropic Hilbert spaces. In \cite{BD3},
we apply our results to construct  transverse slices for the action of the
group of contact diffeomorphisms on the space of compatible Cauchy-Riemann
structures on  three dimensional compact contact manifolds.

While our main interest is in the study of the action of contact
diffeomorphisms on the space of compatible CR~structures, the space of
contact diffeomorphisms is of independent interest. Smooth contact
diffeomorphisms were first studied by Gray \cite{gray:1959a}. Later, Omori
\cite{Omori:1970,Omori:1974} worked within the category of ordinary Sobolev
spaces to show that the space of contact diffeomorphisms is an ILH
(\emph{Inverse Limit Hilbert}) Lie group.
We note that, in the special case where the contact
manifold admits a free transverse $S^{1}$ action, Biquard \cite{Bi} used a
different method to obtain a local parameterization for contact
diffeomorphisms near the identity.

In this paper, we study the group of contact diffeomorphisms from the
perspective of Folland-Stein function spaces, the natural function spaces
respecting the contact structure.  As we show in \cite{BD3}, by working in
Folland-Stein spaces, we are able to obtain normal form theorems for
the action of CR~diffeomorphisms on CR~structures with only
finite regularity.

The paper is structured as follows. In Section~\ref{sec:top-group-str}, we
introduce the anisotropic function spaces, and show that for $s \ge (2n+4)$,
they form an algebra. We believe that the sharp result here would be that
the intersection of $\Wk s$ with the space of bounded functions is an
algebra for all $s$; however, the result which we have stated is sufficient
for our purposes.  In Section~\ref{sec:rumin}, we review Rumin's complex,
and state the results which we will use in the analysis of the space of
contact diffeomorphisms. In Section~\ref{sec:param}, we show that Ebin's
trick of constructing a local coordinate system for the space of
diffeomorphisms works equally well for the anisotropic Folland-Stein spaces,
and then use Rumin's estimates to show the space of contact diffeomorphisms
locally forms a smooth Hilbert submanifold within the coordinate
chart. Using this result, we easily prove that the space of contact
diffeomorphisms is a smooth Hilbert manifold modelled on the anisotropic
function spaces.

\subsection{Notation}
\label{notation}
We summarize here the notation and conventions used throughout the paper.

If $A$ is a subset of a topological space $X$, then $\overline{A}$ denotes
 the closure of $A$ in $X$. If $A$ and $B$ are subsets of $X$, then  the
 notation $A\Subset B$ means that $A$ is compactly contained in
 $B$.

 If $\mathcal{F}:\mathcal{A} \to \mathcal{B}$ is a map between Banach
 spaces, with norms $\|\,\cdot\,\|_{\mathcal{A}}$ and
  $\|\,\cdot\,\|_{\mathcal{B}}$, respectively, then the expression
 \[
		     \|\mathcal{F}(f)\|_{\mathcal{B}} \prec
		     \| f \|_{\mathcal{A}}
 \]
 means that there is a constant $C>0$ such that
 $\|\mathcal{F}(f)\|_{\mathcal{B}} \le C \| f \|_{\mathcal{A}}$ for all $f \in
 \mathcal{A}$.

We give $\R^m$ the standard inner product $\inner{\cdot,\cdot}$, and we let
$\left|\cdot\right|$ denote the corresponding norm. The symbols
$\inner{\cdot,\cdot}_{s}$ and $\norm{\cdot}_{s}$, for $s=0,1,\dots$ denote
the Folland-Stein inner products and norms, respectively.

If $N$ is a smooth manifold, then $TN$ and $T^*N$ denote its tangent and
cotangent bundles, respectively; $\Lambda^p N$ denotes
the $p$-th exterior power of $T^*N$;  $\W^{p}(N)$ denotes the space of
smooth $p$-forms on $N$;  $\Lie_{X}\beta$ denotes the Lie derivative of the
form $\beta$ with respect to the vector field $X$; and $X\inter \beta$
denotes interior evaluation. If $N$ has a Riemannian metric, then $|X|$
denotes the norm of the tangent vector $X$ with respect to that metric.

The symbol $C^r(E)$, $r = 0,1,2,\dots,\infty$, denotes the space of
$C^r$-sections of a fiber bundle $E\to M$, equipped with the topology of
uniform convergence of derivatives up to order $r$ on compact
sets. Similarly, $C^r(M,N)$ denotes the space of $C^r$ maps from $M$ to $N$

We endow $\R^{2n+1}$ with
the contact structure defined by the one-form
\[
	 \n_{0} =  dx^{2n+1} -  \sum_{j=1}^{n} x^{n+j} dx^j \,,
\]
where  $(x^1,\dots,x^n,x^{n+1},\dots,x^{2n},x^{2n+1})$ are the standard
coordinates on $\R^{2n+1}$, and we let $dV_0$ denote the standard volume
form:
\[
        dV_0 = \frac{1}{n!} \n_0\wedge (d\n_0)^n\,.
\]
We  denote the contact distribution of $\n_{0}$ by  $H_0 \subset T\R^{2n+1}$
and we set
\[
  T_{0} = \pd{}{x^{2n+1}}\,, X_j = \pd{}{x^j} + x^{n+j}\, \pd{}{x^{2n+1}}\,, \text{ and } X_{n+j} = \pd{}{x^{n+j}}
\,, 1 \leq j \leq n \,.
\]
Observe that the collection $\{X_j, 1\leq j\leq 2n
\}$ is a global framing for $H_0$.
Notice that the 1-forms
\[
	 \n_0,  dx^{j},\, dx^{n+j},\, 1 \leq j \leq n,
\]
are the dual coframe to $T_0$, $X_j$, $X_{n+j}$, $1 \leq j \leq n$.

Let $f = (f^1,\dots,f^m)$ be a smooth, $\R^m$-valued function defined on the
closure of a domain $D \Subset\R^{2n+1}$. We define
 \[
	X_{I} f = \begin{cases}
		     X_{i_1}X_{i_2}\dots X_{i_t}f &\text{for $t>0$}\\
		     f &\text{for $t=0$} \,,
		  \end{cases}
\]
where we have introduced the multi-index notation $ I=(i_1,\dots,i_t)$,
$1\leq i_j\leq 2n$ and $X_{I}f = (X_{I}f^1,\dots,X_{I}f^m)$. (For $t=0$,
$I$ denotes the \emph{empty index} $I=()$.)  The integer $t$ is called the
\emph{order} of $I$ and written $|I|$.

Throughout this paper, $M$ denotes a fixed smooth, compact contact manifold of
dimension $2n+1$, with contact distribution $H \subset TM$. We call
sections of $H$ \emph{horizontal vector fields}. 
We let
\[
   \piH : T^{*}M \to H^{*} 
\]
denote the projection map; by abuse of
notation, we also let
\[
       \piH: \Lambda^{p} M \to \Lambda^p H^{*}
\]
denote the extension of $\piH$ to the exterior product bundles.
For convenience, we assume that $M$ supports a fixed contact
one-form\footnote{None of our results depend on this assumption, for if the
  line bundle $TM/H$ is non-trivial, we can lift to a double cover of $M$,
  where a global contact form does exist.} $\n$.  The {\em characteristic}
(or {\em Reeb}) vector field $T$ is the unique vector field satisfying the
conditions $T\inter \n = 1$ and $T\inter d\n = 0$. We can then identify the
dual contact distribution with the annihilator of $T$, i.e.
\[
     H^{*}
= \{ \b \in T^{*}M \st   T \inter \b = 0 \} \subset T^{*}M \,;
\]
more generally
\[
          \Lambda^{p}H^{*} = \{ \b\in \Lambda^{p}(M) \st T \inter \b = 0 \}\,,
\]
and we have the identity
\begin{equation}
\label{eq:piH-formula}
           \piH(\b) = T \inter ( \n \wedge \b) \,.
\end{equation}
Two  forms $\b_1$ and $\b_2$ on $M$ are said to be \emph{equal mod $\n$},
written
\[
	\b_1 = \b_2 \modeta \,,
\]
if  and only if $\b_1 = \b_2 + \n \wedge \a$, for some $\a\in \W^*(M)$. An easy
exercise in the exterior calculus proves the equivalences
\[
   \b_1 = \b_2 \modeta 
\iff \n \wedge (\b_1-\b_2) = 0 
\iff  T \inter \left( (\b_1-\b_2)\wedge\n\right) = 0
\iff \piH(\b_1) = \piH(\b_2) \,.
\]
We will also fix an endomorphism $J : H \rightarrow H$ such that $J^2 =
-Id$ and such that the operator $X \mapsto d\n(X,JX)$ is non-negative.
(Such an endomorphism always exists.) We let let $g$ denote the Riemannian
metric defined by the formula
\[
     g( X, Y) = \n (X) \n(Y) + d\n (X, J Y)\,,
\] 
were we have  extended $J$ to a map $J:TM\to TM$ by setting $J(T)=0$.
The endomorphism $J$ and the metric $g$ are said to be {\em adapted} to the
contact structure.  Finally, $ * $ denotes the Hodge star operator
associated to the metric $g$.

We say that a chart $\phi:U \to \R^{2n+1}$ for $M$ is an \emph{adapted
  coordinate chart} if $\n = \phi^*\n_{0}$. It follows that the identities $\phi_* T =
T_0$ and $\phi_* H = H_0$ hold for  $\phi$  adapted.
  An \emph{adapted atlas} for $M$ is a finite,
smooth atlas $\{\phi_{\a}: U_{\a} \to \R^{2n+1}\}$, consisting of adapted
coordinate charts, together with open regions $D_{\a} \Subset
\phi_{\a}(U_{\a})$ such that $\{ W_{\a}= \phi^{-1}_{\a}(D_{\a})\}$ covers
$M$. By compactness of $M$ and Darboux's Theorem for contact structures
\cite[page 362]{Arn:1978}, $M$ has an adapted atlas.  An \emph{adapted
  coordinate chart} for a fibre bundle $\pi: E\to M$ with $m$-dimensional
fibres is a coordinate chart for $E$ of the form
\[
	   \psi: U \to \phi(V) \times \R^{m} \st q 
	   \mapsto (\phi(\pi(q)), \chi(q) )
\]
with $\psi$ surjective and where $\phi:V \to \R^{2n+1}$ is an adapted chart
for $M$.  The chart is said to be \emph{centered at the point $q_0$} if in addition
$\psi(q_0) = (0,0)$. If $\s:V \to U$ is a local section of $E$, the function
$f_{\s}:\phi(V)\to\R^m$ defined by the formula
\[
	  f_{\s} = \chi\circ\s\circ \phi^{-1}
\]
is called the \emph{local representation} of $\s$.

\section{Global analysis on contact manifolds} 
In this section, we develop some of the analytical machinery we need to
study the space of contact diffeomorphisms of $M$. We begin by introducing
the \emph{Folland-Stein} spaces $\Wk{s}$ associated to a compact contact
manifold. We then introduce the notions of \emph{horizontal jet} of a
section of a fibre bundle and define the notion of \emph{contact order} of
a differential operator.  We close this section with a discussion of
\emph{Rumin's Complex} \cite{Rumin:1994} and its associated Hodge theory,
which we need in Section~\ref{sec:param} to construct a local
parameterization of the group of $\Wk{s}$-contact diffeomorphisms of a
compact contact manifold $M$.  Most of the results in this section are
extensions of definitions and theorems in \cite{Palais:1968} and
\cite{Ebin:1970} to the context of contact manifolds.

\subsection{Folland-Stein function spaces}
\label{sec:folland-stein}
\label{subsec:func-space}

Let  $D$ be a bounded domain in $\R^{2n+1}$. 
The {\em Folland-Stein space} $\Wk{s}(D,\R^m)$ is the Hilbert space
completion of the set of smooth, $\R^m$-valued functions on $\overline{D}$
(the closure of $D$) with respect to the
inner product
\[
    (f,g)_{D,s} :=  \sum_{0 \leq |I| \leq s} \int_{D}
       \inner{X_I f\, X_I g } \,dV_{0} \, .
\]
The associated norm is written $\|f\|_{D,s} = \sqrt{(f,f)_{D,s}}$. When no
  confusion is likely to arise, we suppress reference to $D$ and write
  $\|f\|_{s}$; and we set $\Wk{s}(D) := \Wk{s}(D,\R^1)$.

\medskip

\begin{remark}
\label{rem:invariance-1}
Although we used the contact framing $\{X_j \st 1 \leq j \leq 2n\}$ of
$H_{0}$ and the volume form $dV_{0} = \frac{1}{n!}\n_{0}\wedge(d\n_{0})^n$
to define the inner product, an equivalent norm results if the framing is
replaced by any smooth framing of $H_0$ and $dV_{0}$ is replaced by any
smooth volume form on $\bar{D}$.  In particular, suppose that $D'\subset
\R^{2n+1}$ is another bounded domain and $F:\overline{D} \to \overline{D}'$
is a smooth diffeomorphism that restricts to a contact diffeomorphism
between $D$ and $D'$. Then composition with $F$ induces an isomorphism
\[
       \Wk{}_{F}: \begin{cases}\Wk{s}(D',\R^m) &\to \Wk{s}(D,\R^m)\\
                                          f   & \mapsto f \circ F
\end{cases}
\]
between Banach spaces. To see this, notice that because the derivative of
$F$ respects the contact distribution, $X_{I}(f\circ F)$ is a
linear combination of terms of the form $ c_{J} \cdot(X_Jf)\circ F$, $|J|
\leq s$, where
$c_{J}$ denotes a smooth function formed from $F$ and its derivatives. It
follows that $f\circ F$ is of class $\Wk{s}$.
We caution the reader that the condition that $F$ be a contact
diffeomorphism is essential. For if $F$ does not preserve the contact
distribution then the expansion of $X_{I}(f\circ F)$ will in general
involve terms of the form $X_{J}f$ with $|J|>s$, which may not be square
integrable.
\end{remark}

 The next lemma follows immediately from the definition of $\Wk{s}$.
\begin{lemma}
\label{lemma:induction}
The estimate
\[
      \| f \|_{s} \prec \sum_{j} \| X_{j}f \|_{s-1} + \|f\|_{0}
\]
is satisfied for all $f \in \Wk{s}(D)$, $s > 0$.
\end{lemma}

We shall repeatedly make use of the following Sobolev Lemma for
Folland-Stein spaces, which is an immediate corollary of
\cite[Theorem~21.1]{FS}.

\begin{lemma}
\label{lemma:sobolev}
Let $D' \Subset  D \Subset \R^{2n+1}$, and let $s = k + n+2$, $k \geq 0$.
Let $f \in \Wk{s}(D,\R^m)$. Then the functions $X_{I}f$ are continuous on
$\bar{D'}$ for all multi-indices of order $|I| \leq k$. Moreover,
\[
	  \max_{x\in \bar{D'}} |X_{I}f (x) | \prec \|f\|_{D,s}
\]
for all $|I| \leq k$.
If $\|f\|_{D,s}<\infty$ for $s = 2k+ n+2$,  then $f$ is of class $C^{k}$
on $\bar{D'}$. Moreover,  the linear map
\[
	 \Wk{s}(D,\R^m) \to C^{k}(\overline{D}',\R^m)
\]
defined by restriction to $\overline{D}'$ is continuous.
\end{lemma}

\subsection{Estimates for algebraic operations}
\label{sec:alg-estimates}

In this section, we prove some basic estimates.  Lemmas~\ref{lemma:algebra}
and \ref{lemma:division} are needed for our proof in
Section~\ref{sec:top-group-str} that composition and inverses of contact
diffeomorphisms are continuous
operations. Lemma~\ref{lemma:loc-poly-estimate} and
Proposition~\ref{prop:palais-axiom5} are fundamental estimates used
throughout the paper.

\begin{lemma}
\label{lemma:algebra}
Let $s \ge 2n+3, k \le s$ and consider open sets $D' \Subset
D\Subset\R^{2n+1}$.  Then, for any functions $f \in \Wk{s}(D), g \in
\Wk{k}(D)$,
\[
  \|f\cdot g\|_{D',k} \prec   \|f\|_{D,s} \cdot \|g\|_{D,k} \,.
\]
Consequently, multiplication extends to a smooth bilinear mapping
\[
      \Wk{s}(D) \times \Wk{k}(D) \rightarrow \Wk{k}(D')\,.
\]
\end{lemma}

\begin{proof}
We need only prove the estimate for smooth functions $f$ and $g$ on
 $D$. Recall that
\[
   \| f \cdot g \|^2_{D',k}
= \sum_{|I| \leq k } \int_{D'} |X_{I}(f \cdot g)|^2
dV_{0}  \,.
\]
Applying the Leibniz rule, we find that
\[
   \| f \cdot g \|^2_{D',k} \prec \sum_{|J| + |K| \leq k } \int_{D'}
  \left| X_{J}(f) \right|^2 \left| X_{K}(g) \right|^2
dV_{0} \,,
\]

There are two cases to consider: $k \le (n+1)$ and $k>n+1$. In the first case
 $|J| \leq n+1$ for every summand, and we have the estimates
\[
\begin{split}
 \int_{D'}
  \left| X_{J} f \right|^2  \left| X_{K} g  \right|^2 dV_{0}
  \prec& \left(\sup_{x\in D'} | X_{J}f (x)|^2 \right) \,
		      \| X_{K} g  \|_{D',0}^{2}\\
  \prec&  \|f\|_{D,|J|+n+2}^{2} \, \|g\|_{D,k}^{2}\\
  \prec&  \|f\|_{D,s}^{2} \, \|g\|_{D,k}^{2} \,,
\end{split}
\]
where we have used the Sobolev inequality\footnote{We note here that
when $|J|=n+1$, $|J|+n+2 = 2n+3$; whence the condition $s\geq 2n+3$
in the statement of the Lemma.} (Lemma~\ref{lemma:sobolev}) at
the penultimate inequality.

In the latter case, in each term either $|J| \leq n+1$ or $|K| \leq
k-(n+2)$.  In the first instance we bound the term by $\|f\|_{D,s}^{2}
\, \|g\|_{D,k}^{2}$ as before; in the latter case we have
\[
\begin{split}
 \int_{D'}
  \left| X_{J} (f) \right|^2  \left| X_{K} (g)  \right|^2 dV_{0}
  \prec&  \| X_{J} f  \|_{D',0}^{2}\, \sup_{x\in D'} | X_{K}g (x) |^2\\
  \prec&  \|f\|_{D,|J|}^{2} \, \|g\|_{D,|K|+n+2}^{2}\\
  \prec&  \|f\|_{D,s}^{2} \, \|g\|_{D,k}^{2} \,,
\end{split}
\]
where we have again made use of Lemma~\ref{lemma:sobolev}.
Summing over all terms gives the final estimate.
\end{proof}

\begin{lemma}
\label{lemma:division}
Let $D'$ and $D$ be open sets with $D' \Subset D \Subset \R^{2n+1}$, and
let $f$ be a function in $\Wk{s}(D)$, where $s \geq 2n+3$. Suppose that
$1/f$ is bounded from above on $\bar{D}'$ by a positive constant
$C>0$. Then
\[
       \| 1/f \|_{D',s}  \prec  \left(1 + \|f\|_{D,s} \right)^{s}\,.
\]
Consequently, $1/f$ is contained in $\Wk{s}(D')$.
Moreover if $1/f' <C$ for another function $f'\in \Wk{s}(D)$ then
\[
	\| 1/f - 1/f' \|_{D',s} 
      \prec (1 + \|f\|_{D,s})^s (1 + \|f'\|_{D,s})^s \, \|f -  f'\|_{D,s} \,.
\]
\end{lemma}

\begin{proof}
 We have to estimate the quantities
\[
	    \int_{D'} \left| X_J \left(\frac{1}{f} \right) \right|^2
dV_{0}
	    \,,
\]
for $|J| = t \leq s$.  Now, by the quotient and product rules each
such term is bounded by a sum of expressions of the form
\[
     \int_{D'} \left| \frac{ X_{J_{1}}f \cdot X_{J_{2}}f \cdot \dots
     X_{J_{p}}f }{ f^{p+1} } \right|^2
     \, dV_{0} \,,
\]
where $|J_1| + |J_2|+\dots+ |J_p| = t$.  Notice
that $|J_j| > s/2$ for at most one multi-index. Hence, computing as in
Lemma~\ref{lemma:algebra}, we have
\[
\begin{split}
     \int_{D'} \left| \frac{ (X_{J_{1}}f) \cdot (X_{J_{2}}f)  \dots
     (X_{J_{p}})f }{ f^{p+1} } \right|^2 \, dV_{0}
 \prec\; &
       C^{2(p+1)} \cdot \left( \| f\|^{2}_{D,s} \right)^{p-1}
     \|f\|^2_{D,s}\\
 \prec\; & (1 + \|f\|_{D,s} )^{s} 
\end{split}
\]
Summing over all terms gives the first estimate.

The second estimate follows immediately by applying the first estimate and
applying Lemma~\ref{lemma:algebra} to the quantity $1/f - 1/f' = (f'-f)/f f'$.
\end{proof}

A minor modification of the proof of Lemma~\ref{lemma:algebra} gives an
estimate for the product of several functions.
\begin{lemma}
\label{lemma:loc-poly-estimate}
Let $D' \Subset D\Subset\R^{2n+1}$, with $D'$ and $D$ open and
 $s \geq 2n+4$.  Then 
\[
  \|f_1\cdots f_p\|_{D',s} \prec   \norm{f_1}_{D,s} \dots \norm{f_p}_{D,s}
\]
for  all $f_{j} \in \Wk{s}(D)$ for $j=1,2,\dots,p$.

Moreover for $s > 2n+4$,
\[
 \|f_1\cdots f_p\|_{D',s} \prec 
 \sum_{j=1}^p  \norm{f_1}_{D,s-1} \dots
\norm{f_{j-1}}_{D,s-1}\norm{f_{j}}_{D,s} \norm{f_{j+1}}_{D,s-1}\dots
\norm{f_p}_{D,s-1}
\,,
\]
for all $f_j \in \Wk{s}(D)$, $j=1,2,\dots,p$.
\end{lemma}
\begin{proof}
The proof is similar to the proof of Lemma~\ref{lemma:algebra}. By the
product rule,
\[
\norm{f_1\dots f_p}^2_{D',s} 
\prec \sum_{|J|\leq s} \int_{D'} |X_J(f_1\cdots f_p)|^2\, dV_{0}
\prec \sum_{|J_1| + \dots +|J_p| \leq s} \int_{D'} |X_{J_1}f_1|^2 \cdots
|X_{J_p}f_p|^2 \,dV_{0}  \,,
\]
where $X_I f$ is defined in Section~\ref{sec:folland-stein}.
We need only bound each term in
the right-hand summation. 
Since $s \geq 2n+4$, it follows that $|J_j| > s/2 \geq n+2$ for at most
one multi-index, say $J_j$ 
in the right-hand sum and that (since $n+2 \leq s/2$)
\[
   |J_i|+n+2 \leq s/2 + n+2 \leq s \text{ for } i \neq j \,.
\] 
Hence by  Lemma~\ref{lemma:sobolev}, $X_{J_i}f_i$ is continuous and
	   $\sup_{x\in D'} |X_{J_i}f_i| \prec \norm{f_i}_{D,s}$.
Consequently,
\[
 \int_{D'} |X_{J_1}f_1|^2 \dots |X_{J_p}f_p|^2 \,dV_{0} 
\prec \prod_{i\neq j} \sup_{x\in D'} |X_{J_i}f_i|^2 \,\cdot \int_{D'}
|X_{J_j}f_j|^2\,dV_{0}
\prec \norm{f_1}_{D,s}^2 \dots \norm{f_{p}}_{D,s}^2 \,,
\]
from which the first estimate follows.

Now suppose that $s > 2n+4$. Then in the previous paragraph
 $|J_i| + n + 2 \leq s-1$; and
   $\sup_{x\in D'} |X_{J_i}f_i| \prec \norm{f_i}_{D,s-1}$ for $i\neq j$,
 yielding the estimate
\[
\begin{split}
 \|f_1\cdots f_p\|^2_{D',s}
 \prec 
 \sum_{j=1}^p  \norm{f_1}^2_{D,s-1} \dots
\norm{f_{j-1}}^2_{D,s-1}\norm{f_{j}}^2_{D,s} \norm{f_{j+1}}^2_{D,s-1}\dots
\norm{f_p}^2_{D,s-1}\\
 \prec \left( \sum_{j=1}^p  \norm{f_1}_{D,s-1} \dots
\norm{f_{j-1}}_{D,s-1}\norm{f_{j}}_{D,s} \norm{f_{j+1}}_{D,s-1}\dots
\norm{f_p}_{D,s-1} \right)^2\,.
\end{split}
\] 
\end{proof}

\subsection{The Folland-Stein space of sections of a  vector bundle}
\label{subsec:vector-bundle}
In this section we define the Folland-Stein of sections of a vector
bundle over a contact manifold.

We begin by extending the  definition of the Folland-Stein space
$\Wk{s}(D,\R^m)$  of functions to the space
$\Wk{s}(M,\R^m)$ of functions on a compact contact manifold.   Let
$(\phi_{\a}, U_{\a}, D_{\a})$ be an adapted atlas for $M$.  A function
$f:M \to\R^m$ is said to be a {\em $\Wk{s}$-function} if the functions
$f_{\a} = f \circ \phi_{\a}^{-1}$ lie in $\Wk{s}(D_{\a},\R^m)$ for all
$\a$. The formula
\[
      \inner{f,g}_{M,s} = \sum_{\a} \inner{f_\a,g_\a}_{D_{\a},s}
\]
makes $\Wk{s}(M,\R^m)$ into a separable Hilbert space.
 The Sobolev
Lemma~\ref{lemma:sobolev} clearly extends to this setting:

\begin{lemma}
\label{lemma:sobolev-global}
Let $s = k + n+2$, and let $Y_j$, $j=1,2,\dots k$ be smooth sections of
$H\subset TM$.  Then for any function $f \in \Wk{s}(M,\R^m)$, the functions
$Y_1 Y_2\dots Y_kf$ are continuous on $M$. Moreover,
\[
	  \max_{x\in M} |Y_1 Y_2 \dots Y_k f (x) | \prec \|f\|_{M,s}
\]
If $\|f\|_{M,s}<\infty$ for $s = 2k+ n+2$,  then $f$ is of class $C^{k}$
on $M$. In particular for $s=2k+n+2$,  the linear map
\[
	 \Wk{s}(M,\R^m) \to C^{k}(M,\R^m)
\]
is continuous.
\end{lemma}
Similarly, Lemma~\ref{lemma:loc-poly-estimate} assumes the following global
form:
\begin{lemma}
\label{lemma:global-poly-estimate}
If $s \geq 2n+4$
\[
  \|f_1\cdots f_p\|_{M,s} \prec 
  \norm{f_1}_{M,s} \dots \norm{f_p}_{M,s} \,,
\]
and if $s > 2n+4$,
\[
 \|f_1\cdots f_p\|_{M,s} \prec 
 \sum_{j=1}^k  \norm{f_1}_{M,s-1} \dots
\norm{f_{j-1}}_{M,s-1}\norm{f_{j}}_{M,s} \norm{f_{j+1}}_{M,s-1}\dots
\norm{f_p}_{M,s-1}
\,,
\]
for all $f_j \in \Wk{s}(M)$, $j=1,2,\dots,p$.
\end{lemma}

\smallskip

We define the Hilbert space of \emph{Folland-Stein sections} $\Wk{s}(E)$
for $\pi:E\to M$ a smooth vector bundle of rank $m$ as follows.  View
$\Wk{s}(M,\R^r)$, $r\geq 1$, as the Folland-Stein space of sections of the
trivial vector bundle $M\times\R^r\to M$.
For $r$ sufficiently large, there is a vector bundle injection
$\displaystyle E \stackrel{\iota}{\hookrightarrow} M \times \R^r$. Define
an inner product on $C^{\infty}(E)$ by the formula
\[
    \inner{f,g}_{E,s} = \inner{f\circ \iota, g\circ\iota}_{M,s}
\]
for $f,g \in C^{\infty}(E)$, and let $\Wk{s}(E)$ be the Hilbert space
completion of $C^{\infty}(E)$ with respect to this inner product.
It is not difficult to check that, although the inner product depends on
$\iota$, the space $\Wk{s}(E)$ does not. 

\begin{remark}
\label{rem:sobolev-bundle-version}
Because $\Wk{s}(E)$ is a closed subset of $\Wk{s}(M,\R^r)$, the Sobolev
Lemma \eqref{lemma:sobolev-global} extends to this setting,
\end{remark}

The next proposition is the analogue of ``Axiom B2'' of Palais (see
\cite[page 10]{Palais:1968}) in the setting of contact manifolds.

\begin{proposition} 
\label{prop:palais-axiom2}
Let $F:M \to N$ be a smooth contact
  diffeomorphism between two compact, contact manifolds. Let $E\to N$ be a
  smooth vector bundle over $N$ and let $F^*E \to M$ be its pull-back to
  $M$. Then  the map $\s \mapsto \s \circ F$ is a
  Hilbert space isomorphism between $\Wk{s}(E)$ and $\Wk{s}(F^*E)$.
\end{proposition}

\begin{proof} 
Choose a bundle injection $\iota:E\hookrightarrow N\times\R^r$, and
let  $\{\phi_{\a}: U_{\a} \to \R^{2n+1}\}$ be an adapted atlas for $N$ (see
Introduction). Set $U'_{\a} = F^{-1}(U_{\a})$ and $\phi'_{\a} = \phi_{\a}
\circ F$. Then  $\{\phi'_{\a}: U'_{\a} \to \R^{2n+1}\}$ is an adapted atlas
for $M$. Use this atlas to define the inner product on $\Wk{s}(M,\R^r)$.
Then by construction
\[
         (f\circ F, g\circ F)_{s,M} =   (f, g)_{s,N}  
\]
for all $f,g \in \Wk{s}(N,\R^r)$.  Restricting $f$ and $g$ to sections of
$E$ then gives the result.
\end{proof}

The next proposition shows that $\Wk{s}$ satisfies ``Axiom~B5'' of Palais
(\cite[page 39]{Palais:1968} for all $s \geq 2n+4$.

\begin{proposition}
\label{prop:palais-axiom5}
Let $E_j\to M$, $j=1,2$ be smooth vector bundles over $M$ and let $F:
E_1\to E_2$ be a smooth (not necessarily linear) fibre-preserving map.
Then the map
\[
    \Wk{}_{F}: \Wk{s}(E_1) \to \Wk{s}(E_2) \,:\,
		 \s \mapsto  F \circ \s
\]
is a $C^{\infty}$ map for all $s \geq 2n+4$.

Moreover, if $E_1$ is equipped with norm $|\cdot|_{E_1}$, then
for every $\s \in \Wk{s}(E_1)$ and every $c > 0$ there is a polynomial
 $Q$  with non-negative coefficients of total degree at most
 $s$ such that
\begin{equation}
\label{poly-estimate}
\tag{*}
		\|F\circ \s - F\circ \s'\|_{s} \le \|\s - \s' \|_{s} Q^{s}
                (\|\s \|_{s} , \|\s' \|_{s})
\end{equation}
for all $\s' \in \Wk{s}(E_1)$, with $|\s' - \s| <c$.
\end{proposition}

\begin{proof}
 Our proof  follows a similar argument  in \cite[Theorem~11.3]{Palais:1968}.
 We first show that $\Wk{}_{F}$ satisfies the polynomial estimate
 \eqref{poly-estimate}, from which continuity of $\Wk{}_{F}$ follows.  It
 suffices to work in local coordinates of an adapted atlas for $M$. Then
 $\Wk{s}$ sections of $E_j$ can be identified with elements of
 $\Wk{s}(D,\R^{m_j})$, where $m_j$ is the fibre dimension of $E_j$.  Let
 $f:D \to \R^{m_1}$ be the local coordinate representation of $\s$, and
 choose a constant $c>0$. Choose $\s'$ so that it's local representative $g
 \in \Wk{s}(D,\R^{m_1})$ satisfies $\sup_{x \in D} |g(x)-f(x)| < c$.  Then
 in local coordinates $F\circ\s(x) = F(x,f(x))$ and $F\circ\s'(x) =
 F(x,g(x))$. By smoothness of $F$ and compactness of $\overline{D'}$, all
 derivatives of $F$ are bounded and smooth on the set $\{(x,y) \st x \in
 \overline{D'}, |y - f(x)|<c$.  Hence, there is a fixed constant $C>0$ such
 that
\begin{equation}
\tag{i}
         |F^{(k)}(x,y)| < C \text{ and }  | F^{(k)}(x,y) - F^{(k)}(x,z) | < C |y-z|
\end{equation}
for all $(x,y)$ and $(x,z)$ with $|y-f(x)| \leq c$ and $|z-g(x)| \leq
c$, where $F^{(k)}(x,y)$ denotes any mixed partial
derivative of $F$ of order $k\leq s$.

 Next recall that
 $\norm{F\circ f - F\circ g}^2_{D',s}$  is a sum of integrals of the form
\begin{equation}
\tag{ii}      \int_{D'} |X_I  \left\{F(x, f(x)) - F(x,g(x)) \right\}|^2 \, d V
\end{equation}
for $|I|\leq s$.  By the chain rule, $X_I \left\{F(x,f(x)) - F(x,g(x))
\right\}$ is a finite sum of terms of the form
\begin{multline}
\tag{iii}   
 F^{(k)}(x,f(x)) \cdot X_{I_1}f(x) \cdots X_{I_k}f(x) -
 F^{(k)}(x,g(x)) \cdot X_{I_1}g(x) \cdots X_{I_k}g(x)
\\
 =\left\{F^{(k)}(x,f(x)) - F^{(k)}(x,g(x)) \right\}\cdot X_{I_1}f(x) \cdots X_{I_k}f(x)\\
 +
 (F^{(k)}(x,g(x)) \cdot \left\{ X_{I_1}f(x) \cdots X_{I_k}f(x) -  X_{I_1}g(x)
 \cdots X_{I_k}g(x) \right\}
\end{multline}
where $0 \leq k \leq s$ and $\sum_{j=1}^{k} |I_j| \leq s$.  Applying (i) to
(iii) gives the
estimate
\begin{multline}
\tag{iv}
| F^{(k)}(x,f(x)) \cdot X_{I_1}f(x) \cdots X_{I_k}f(x) -
 F^{(k)}(x,g(x)) \cdot X_{I_1}g(x) \cdots X_{I_k}g(x)| \\
\leq C\, |f(x)-g(x)| \cdot |X_{I_1}f(x) \cdots X_{I_k}f(x)| + 
     C \,| X_{I_1}f(x) \cdots X_{I_k}f(x) -  X_{I_1}g(x)
 \cdots X_{I_k}g(x)|
\end{multline}
The right-hand side of (iv) can, in turn, be bounded by a finite sum of terms
form
\begin{equation}
\tag{v}
        C |X_{I_0}( f(x)-g(x))| \cdot |X_{I_1}f(x)| \dots |X_{I_{k'}}f(x) |
\cdot  |X_{I_{k'+1}}g(x)| \dots |X_{I_{k'+k''}}g(x)| 
\end{equation}
where $0 \leq |I_i|$ and $\sum_{i=0}^{k'+k''} |I_i| \leq s$.
Substituting (v) into (ii) shows that $\norm{F\circ f - F\circ g}^2_{D',s}$
is bounded by a sum of integrals of the form
\begin{multline}
\tag{vi}
 \int_{D} C^2  |X_{I_0}( f(x)-g(x))| \cdot |X_{I_1}f(x)| \dots |X_{I_{k'}}f(x) |
\cdot  |X_{I_{k'+1}}g(x)| \dots |X_{I_{k'+k''}}g(x)|\\
 |X_{J_0}( f(x)-g(x))| \cdot |X_{J_1}f(x)| \dots |X_{J_{\ell'}}f(x) |
\cdot  |X_{J_{\ell'+1}}g(x)| \dots |X_{J_{\ell'+\ell''}}g(x)|  \, dV_{0}
\end{multline}
where $0 \leq |J_j|$ and $\sum_{j=0}^{\ell'+\ell''} |J_j| \leq s$.

 Notice that
   $|I_i| \geq s/2$ and $|J_j| \geq s/2$ for at most one $i$ and at most one $j$,
  and since $s \geq 2n+4$, the Sobolev Lemma~\ref{lemma:sobolev} applies, to show that
the remaining factors in the integrand are all continuous, hence bounded on
the compact set $\{(x,y) \st x \in \overline{D'}, |y-f(x)| \leq c \}$. It
follows that the integral in (vi) is bounded by an expression of the form
\[
         C' \norm{f-g}^2_{D,s} \norm{f}_{D,s}^{k'} \norm{g}_{D,s}^{k''}
           \norm{f}_{D,s}^{\ell'} \norm{g}_{D,s}^{\ell''} \,,
\]
for $C'$ a constant depending on $s$, $F$, $c$, and $D'$. Since $k'+k''
\leq s$ and $\ell' + \ell'' \leq s$, it follows that
\[
               \norm{F\circ f - F\circ g}_{D',s} \prec
                   \norm{f-g}_{D,s} \; Q(\norm{f}_{D,s}, \norm{g}_{D,s} )\,,
\]
where $Q(u,v)$ is a polynomial of bidegree at most $s$ in
$u$ and $v$, with non-negative coefficients. Applying this to each chart in
an adapted atlas yields the global estimate \eqref{poly-estimate}.

We now show that $\Wk{}_{F}$ is $C^1$ with derivative given
by the formula
\[
      d\Wk{}_F = \Wk{}_{\delta F}: \Wk{s}(E_1) \times \Wk{s}(E_1)
\to \Wk{s}(E_2)
\]
where $\delta F : E_1 \times_{M} E_1 \to E_2$ is the smooth fibre bundle map
defined by the formula
\[
       \delta F_x(u,v) = \left. \frac{d}{dh}\right|_{h=0}  F(u + h v)\,,
\text{ for all $x\in M$ and $u,v \in E_{1,x}$.}
\]
To show that
\[
\label{eq:p1}
      \lim_{v\to 0} \frac{
	 \| \Wk{}_{F}(\s +  v) - \Wk{}_{F}(\s) - \Wk{}_{\delta F}(\s,v)
       \|_{s} }{\|v\|_{s}} = 0 \, ,
\]
first observe that $\delta F$ can be expressed as a smooth
map of the form
\[
       \delta F : E_{1} \to \text{Hom}(E_1,E_2) \,.
\]
Hence,
\[
     \Wk{}(\delta F): \Wk{s}(E_1) \to \Wk{s} \left(\text{Hom}(E_1,E_2) \right)
\]
is continuous, and for all $\epsilon>0$ there is a $\delta_1$ such that
\[
 \| \Wk{}(\delta F)(\s + v) - \Wk{}(\delta F)(\s) \|_{s} < \epsilon
\text{ whenever $\|v\|_{s} < \delta_1$.}
\]
Using this observation, we compute as follows:
\begin{multline*}
 \Wk{}_{F}(\s +  v) - \Wk{}_{F}(\s) - \Wk{}_{\delta F}(\s,v)
= F(\s + v) - F(\s) - \delta F_{\s} \cdot v \\
= \int_{0}^{1} \left\{ \delta F_{\s + t v} \cdot v - \delta F_{\s} \cdot v
\right\} \,dt
=
 \int_{0}^{1} \left\{ \Wk{}_{\delta F} (\s + t v) (v) - \Wk{}_{\delta
 F}(\s) ( v )
\right\} \,dt 
\end{multline*}
Hence, if $\|v\|_{s} < \delta_1$ then
\[
  \| \Wk{}_{F}(\s +  v) - \Wk{}_{F}(\s) - \Wk{}_{\delta F}(\s,v) \|_{s}
\leq\;
   \int_{0}^{1} \| \Wk{}_{\delta F} (\s + t v) (v) - \Wk{}_{\delta F}(\s) (
v ) \|_{s} \,dt
<\;
     \epsilon \|v\|_{s}\,.
\]
This show that $\Wk{}_{F}$ is  differentiable at
$\s$.  That it is continuously differentiable follows from the identity
$d\Wk{}_{F} = \Wk{}_{\delta F}$ and continuity of $\Wk{}_{\delta F}$.

That $\Wk{}_{F}$ is smooth follows by induction.
For assume  that for some $k>0$, $\Wk{}_{F}$ is $C^{k}$, for
all smooth $F:E_1\to E_2$, and all $E_j$.
To show that $\Wk{}_{F}$ is $C^{k+1}$, we need only show that its derivative
$d\Wk{}_{F}: \Wk{s}(E_1\times E_1) \to \Wk{s}(E_2)$ is $C^k$. But
$d\Wk{}_{F} = \Wk{}_{\delta F}$, and $\delta F$ is, a smooth fibre bundle
map. Consequently, $d\Wk{}_{F}$ is $C^k$, completing the induction step.
\end{proof}

\subsection{The Folland-Stein space of sections of a  fibre bundle}
\label{subsec:fibre-bundle}

In this section we define  the Folland-Stein space
$\Wk{s}(E)$ of sections of $E$ for $s \geq 2n+4 $ in the case where
 $\pi: E
\to M$ is a smooth fibre bundle over $M$. For this range of $s$,
$\Wk{s}(E)$ is a smooth infinite dimensional manifold modelled on
the Folland-Stein space of sections of  certain vector bundles.
The construction and the proof are due to Palais (see \cite[Chapters 12 and
13]{Palais:1968} for details). 
We emphasize that this construction depends heavily on 
Propositions~\ref{prop:palais-axiom2} and \ref{prop:palais-axiom5}
(Palais' Axioms B2 and B5), which are satisfied for $s \geq 2n+4$.

Let $\pi:E\to M$ be a fibre bundle, and choose 
a smooth section $\s \in C^{\infty}(E)$. 
The \emph{bundle of vertical tangent vectors} along $\s$ is
the vector bundle defined by
\[
     T_{\s}(E)  = \{ X\in TE_{\s(x)} \st x\in M, d\pi (X) = 0\} \,. 
\]
Palais shows that there is a  smooth fibre bundle isomorphism
\begin{equation}
\label{eq:sections}
 \psi_{\s}: T_{\s}(E) \stackrel{\simeq}{\longrightarrow}  
    \mathcal{O}_{\s} \stackrel{\text{open}}{\subset} E  
\end{equation}
where $\mathcal{O}_{\s}$ is a neighbourhood of the image of $\s$. Palais
also shows that the the image of every continuous section of $E$ is
contained in a set of the form $\mathcal{O}_{\s}$ for some smooth section
$\s$.

Consequently every continuous section of $E$ can be identified with a
continuous section of $T_{\s}(E)$ for some $\s\in C^{\infty}(E)$, and
$C^{0}(E)$ can we written as the following union of open sets:
\[
	   C^{0}(E) = \bigcup_{\s\in C^{\infty}(E)} C^{0}(T_{\s}(E)) \,.
\]
Since $s \geq 2n+4\geq n+2$, Lemma~\ref{lemma:sobolev-global} applies to
give continuous inclusions $\Wk{s}(T_{\s}(E)) \subset C^{0}(T_{\s}(E)$.

We may
thus define the \emph{Folland-Stein} space $\Wk{s}(E)$ to be the union
\[
	   \Wk{s}(E) = \bigcup_{\s\in C^{\infty}(E)} \Wk{s}(T_{\s}(E)) \,,
\]
equipped with the weakest topology such that $\Wk{}_{\psi_{\s}}:
\Wk{s}(T_{\s}(E)) \to \Wk{s}(E)$ is a continuous open map for all $\s\in
C^{\infty}(E)$.

In fact, $\Wk{s}(E)$ is a Hilbert manifold with $C^{\infty}$ atlas given by
the charts $\Wk{}(\psi^{-1}_{\s})$; and as Palais shows, smoothness of the
transition functions for this atlas follows from
Proposition~\ref{prop:palais-axiom5}.

The above construction is functorial:

\begin{proposition}[Palais, Theorem 13.4]
\label{prop:palais}
Let $E_j\to M$, $j=1,2$ be smooth fibre bundles over $M$ and let
$F: E_1\to E_2$ be a smooth fibre-preserving map.  Then the map
\[
    \Wk{}(F): \Wk{s}(E_1) \to \Wk{s}(E_2) \,:\,
		 \s \mapsto  F \circ \s
\]
is a $C^{\infty}$ map of Hilbert manifolds for all $s \geq 2n+4$.

\end{proposition}

\begin{remark}
\label{rem:sobolev-fibre-bundle}
  By construction, the Sobolev Lemma~\ref{lemma:sobolev-global} extends to
  define a continuous injection $\Wk{s}(E) \subset C^{k}(E)$, for $s =
  \max(2k+n+2, 2n+4)$ .
\end{remark}

\medskip

The case where $E$ is the trivial fibre bundle $E=M\times N \to M$ is an
important special case:

\begin{definition}
\label{def:space-of-maps}
  Let $N$ be a smooth manifold without boundary. For $s \geq 2n+4$, the
  Folland-Stein space $\Wk{s}(M,N)$ of maps from the contact manifold $M$
  to the manifold $N$ is the Folland-Stein space of sections of the trivial
  fibre bundle $M\times N\to M$.
\end{definition}

\begin{corollary}
\label{cor:left-composition}
Let $F:N \to N'$ be a $C^{\infty}$ map between $C^{\infty}$ manifolds. 
Left composition with $F$ defines a $C^{\infty}$
\[
          L^{s}_F: \Wk{s}(M,N) \to \Wk{s}(M,N') \,:\, G \mapsto F \circ G
\]
for all $s \geq 2n+4$. If $F:N \to N$ is a diffeomorphism of $N$ then
$L_{F}^{s}$ is a diffeomorphism of $\Wk{s}(M,N)$.
\end{corollary}

\begin{proof} View composition with $F$ as a smooth bundle map $(x,y)
  \mapsto (x, F(y))$, and apply Proposition~\ref{prop:palais}. If $F$ is a
  diffeomorphism then $R^{s}_{F^{-1}}$ is the smooth inverse of $R_{F}^{s}$.
\end{proof}


\subsection{Horizontal jets and differential operators}
\label{sec:diff-operators}

The goal of this section is to extend the framework of \cite[Chapter
15]{Palais:1968} to the context of differential operators on contact
manifolds.

  Let $\pi:E \to M$ be a smooth fibre bundle over $M$, and  choose an
adapted coordinate chart 
\[
     \psi:U \to \phi(V)\times \R^m \st p \mapsto (x,y) \,,
\]
as defined in Section~\ref{notation}.
Two smooth local sections
$\s_{i}$, $i=1,2$, of $E$ defined on $V$, with local representations
 $f_i: \phi(V)\to\R^m$,  are said to be
\emph{contact equivalent up to order $k$} at a point $p\in V$
 if and only if
\begin{equation}
\label{eq:jet-equiv}
	 X_{I}f_{1}( \phi(p)) = X_{I}f_{2}(\phi(p))
\end{equation}
for every multi-index $I$ with $0 \leq |I| \leq k$. It is easy to check that
contact equivalence is an equivalence relation and that it is independent
of coordinates.  
The \emph{ horizontal $k$-jet of $\s$ at $p$}, written $j^{k}_{H}\s (p)$,
is the equivalence class of the local section $\s$ at $p\in M$,
 $J^{k}_{H}E$ denotes the space of all horizontal
$k$-jets. 
The map $\pi: J^{k}_{H}E \to M$ defined by
\[
       \pi\left( j^{k}_{H}\s (p) \right) = p 
\]
makes the space of horizontal $k$-jets into a fibre bundle with fibres of
dimension $m\cdot N_k$, where $N_k$ is the number of indices $A$ with
$|A|\leq k$.

\begin{remark}
\label{rem:A-notation}
By virtue of the commutation relations among the vector fields
$T_0$, $X_{1}$,\dots $X_{2n}$,
\emph{any} differential operator of the form
$\displaystyle
		Y_1 Y_2 \dots Y_r
$,
where $Y_1,\dots,Y_r$ are arbitrary vector fields on  an open set $V \subset M$,
can be expressed uniquely in local coordinates as a linear combination of operators of the form
\[
   D_{A} :=  \underbrace{X_1\cdots X_1}_{a_1}, \underbrace{X_2\cdots X_2}_{a_2},\dots,
 \underbrace{X_{2n}\cdots X_{2n}}_{a_{2n}} \underbrace{T_0\cdots T_0}_{a_{2n+1}}
\]
where $A = (a_1,a_2,\dots,a_{2n},a_{2n+1})$, $0 \leq a_j$. Moreover, if
$Y_1,\dots,Y_r$ are all horizontal, then $a_1+a_2+\dots+a_{2n}+ 2 a_{2n+1}
= r$.  The integer $|A| =a_1 + \dots + a_{2n} + 2 a_{2n+1} $ is called the
{\em contact order} of $D_{A}$.  It follows that  $\s_1$ and
$\s_2$ are contact equivalent up to order $k$ at $p$ if and only if
\[
	   D_{A} f_1( \phi(p)) = D_{A} f_2( \phi(p))
\]
for all multi-indices $A$ with $|A|\leq k$.
\end{remark}

\begin{lemma}
\label{lemma:k-jets}
Let $E\to M$ be a smooth fibre bundle with fibre dimension $m$.
Then $J^{k}_{H}E \to M$ is a
smooth fibre bundle. 
 Moreover, if $\s$ is a smooth section of $E$ then
$j^{k}_{H}\s: p \mapsto j^{k}_{H}\s(p)$
is a smooth section of $J^{k}_{H}E$.
If $F:E_1\to E_2$ is a smooth map of fibre bundles,
then so is the map
\[
	 J_H^k(F): J^{k}_{H} E_1 \to  J^{k}_{H} E_2
\st j^k_H \s(p) \mapsto j^k_H (F\circ\s) (p)\,.
\]
This construction is functorial, i.e. $J_H^k(G\circ F) =J_H^k(G)\circ
J_H^k(F)$, for $G:E_2\to E_3$ a smooth map of fibre bundles.
\end{lemma}

\begin{proof}
We give an outline of the proof, leaving some details to the reader.
To define a coordinate chart for $J^k_{H}E$, choose a point $q_0\in E$
and let $p_0 = \pi(q_0)$. Let
\[
   \psi: U \to \phi(V) \times \R^n \st q \mapsto (x,y)
\]
be adapted coordinates for $E$ centered at $q_0$ (see Section~\ref{notation});
and let $\widetilde{U}$ be the set of horizontal $k$-jets of sections
of $E$ defined on $V$ and with values in $U\subset E$.
Viewing $\R^{m\cdot N_k}$ as the space of $m\times N_k$ matrices,
we can define local fibre bundle coordinates
\[
    \widetilde{\psi}: \widetilde{U} \to 
\R^{2n+1} \times \R^{m\cdot N_k } \,
\]
by the formula
\[
      \widetilde{\psi}( j_H^k\s(p)) =    \left((x(p)),(D_{A}f_{\s} (p) ) \right) \,,
\]
where $f_{\s}$ is the local representation of $\s$
(Section~\ref{notation}), $D_{A}f_{\s}$ is the derivative of $f_{\s}$ with
respect to the multi-index $A$, 
 and we have
given the set $\{A \st |A|\leq k\}$ the lexicographical ordering.  We need only
show that $\widetilde{\psi}$ is a bijection between $\widetilde{U}$ and
$\phi(V)\times \R^{m\cdot N_k}$ for $U$ a sufficiently small neighbourhood
of $q_0$.  By the discussion in Remark~\ref{rem:A-notation},
$\widetilde{\psi}$ is injective.

To prove surjectivity, choose a point $S = \{S_{A} \in \R^{m} \st |A|\leq k
\} \in \R^{m\cdot N_k}$, and consider the polynomial section
\[
   y = \s_{S}(x) =
    \sum_{A} \frac{1}{A!}
	S_{A} \cdot (x^1)^{a_1}\dots (x^{2n})^{a_{2n}} (x^{2n+1})^{a_{2n+1}}
 \,.\]
where we have adopted the notation
$A! = a_1! \dots a_{2n}! a_{2n+1} !$.
Now for any point $p = (x^1,\dots,x^{2n+1})$, the map
\[
     S \mapsto \widetilde{\psi}(\s_{S}(p)) = (x, D_{A}\s_S (x))
\]
may be viewed as a linear map $L_{x}:\R^{m \cdot N_k}
\to  \R^{m \cdot N_k}$. Notice that 
$L_{x}$ depends smoothly on $x$.
To show that $L_{x}$ is a bijection for $x$ sufficiently near $0$, we need only
verify that $L_{0}$ is injective. But a straightforward computation shows that
\[
	D_{A'} \s_S (0) = \begin{cases}   S_{A} & \text{ for $A'= A $}\\
					      0    & \text{ otherwise,}
			   \end{cases}
\]
i.e. $L_{0}$ is the identity map.  This shows that $\widetilde{\psi}$
is surjective on the fibre of $J^k_{H}E$ over all points $p$
sufficiently near $p_0$. Thus, after a possible shrinking of $U$,
the map
\[
     \widetilde{\psi}: \widetilde{U} \to \R^{2n+1}  \times \R^{m \cdot N_k}
\]
is a bijection between $\widetilde{U}$ and $\phi(U) \times \R^{m \cdot N_k}$.

By letting $q_0$ vary over all $E$, we obtain a smooth atlas for
$J^k_{H}E$.  We topologize $J^k_{H}E$ by requiring each of the charts
$\widetilde{\psi}$ to be a homeomorphism.  That these charts form a smooth
atlas making $J^{k}_{H}E$ into a smooth manifold follows from the
observation that if $\widetilde{\psi}$ and $\widetilde{\psi}'$ are two
charts then $\widetilde{\psi}' \circ \widetilde{\psi}^{-1}$ is a
diffeomorphism between $\widetilde{\phi}( \widetilde{U} \cap
\widetilde{U}')$ and $\widetilde{\phi}'( \widetilde{U} \cap
\widetilde{U}')$.  The proof follows by standard arguments in advanced
calculus (i.e. the Inverse Function theorem and the Chain rule) and is left
to the reader.

That $J^k_H(F)$ and $j^k_H\s$ are smooth  follows from the Chain Rule, as
does the identity  $J_H^k(G\circ F) =J_H^k(G)\circ J_H^k(F)$.
\end{proof}

\begin{remark} In the special case where $E \to M$ is a smooth vector
  bundle, then so is $J^{k}_{H}E \to M$, with linear structure induced by
  the formula
\[
    a_1 j_H^k\s_1 (p) + a_2 j_H^k\s_2(p) =  j_H^k( a_1 \s_1 + a_2 \s_2 ) (p) \,.
\]
\end{remark}

\begin{lemma}
\label{prop:jet-map}
Let $E\to M$  be a smooth fibre bundle over $M$.
 Then the map
\[
    \Wk{}(\jH{}): \Wk{s}(E) \to \Wk{s-k}(J^{k}_{H}E) \,:\,
		 \s \mapsto  \jH{k}\s
\]
is a smooth map of Hilbert manifolds for all $k$ and $s$ such that $s \geq
2n + 4 + k$.  In the special case where $E$ is a vector bundle,
$\Wk{}(\jH{})$ is a bounded linear map of Hilbert spaces.
\end{lemma}

\begin{proof}
From the discussion in Section~\ref{sec:folland-stein}, it suffices to
analyze $\Wk{}(\jH{})$ in the neighbourhood of a fixed section of $E$. By
definition of the Hilbert manifold structure of $\Wk{s}(E)$ we may, without
loss of generality assume that $E$ is a vector bundle, and we must only show
that the linear map $\Wk{}(\jH{})$ is bounded.  By definition of the inner
product on $\Wk{s}(E)$, the result then follows from a local
computation: Choose an open set $V \Subset M$ and a trivialization
$E_{|\overline{V}} \simeq \overline{V} \times \R^m$.  Then a section $\s$ of $E$ over
$V$ is given by an $\R^m$-valued function $f_{\s}$ and its
horizontal $k$-jet $j_H^k\s:V \to J^{k}_{H}M$ can be identified with the $m
\times N_k$ matrix-valued function
\[
\jH{k}\s = 
  \left( D_{A} f_{\s}\right)_{ |A| \leq N_k}
\]
So
\begin{align*}
    \|\jH{k}\s\|^2_{V,s-k}
&=  \sum_j \sum_{0 \leq |J| \leq s-k} \int_{V}  |X_J (\jH{k}\s^j)|^2 dV_{0}
= \sum_{0 \leq |J| \leq s-k} \sum_{|A| \leq k} \int \left| X_J D_{A} f_{\s}\right|^2 dV_{0} \\
&\leq  \sum_{0 \leq |I| \leq s} \int \left| X_I f_{\s}\right|^2 dV_{0}
\leq  \|\s\|_{V,s}^2 \, .
\end{align*}
\end{proof}

\begin{definition}
\label{def:contact-order}
Let $E_1$ and $E_2$ be smooth fibre bundles over $M$. A {\em
 differential operator of contact order $k$} from $E_1$ to $E_2$ is a
 map of the form
\[
       D: C^{\infty}(E_1) \overset{\jH{k}}{\longrightarrow} C^{\infty}(J^{k}_{H}E_1)
		      \overset{F_{*}}{\longrightarrow} C^{\infty}(E_2) \,.
\]
where $F:J^{k}_{H}E_1 \to E_2$ is a smooth fibre bundle map and $F_{*}$ is
defined by the formula
\[
	     F_{*}(\hat{\s}) = F\circ \hat{\s} \text{ for $\hat{\s}\in
  C^{\infty}(J^{k}_{H}E_1)$}\,.
\]
\end{definition}

\begin{proposition}
\label{diffop-est}
Let $D$ be a differential operator of contact order $k$ as above.
Then $D$ extends to a smooth map
\[
       D: \Wk{s}(E_1) \overset{\jH{k}}{\longrightarrow} \Wk{s-k}(J^{k}_{H}E_1)
		      \overset{F_{*}}{\longrightarrow} \Wk{s-k}(E_2) \,,
\]
for all  $s \geq 2n + 4 + k$. 

In the special case where $E_1$ and $E_2$ are normed vector bundles, for
every section $\s \in \Wk{s}(E_1)$ and every constant $c > 0$, there is a
polynomial $Q$ of degree at most ${s-k}$ such that the estimate
\[
\| D\s - D\s'\|_{s-k} \le 
           \|\s-\s'\|_{s} \cdot Q(\|\s\|_{s}, \|\s'\|_{s})
\]
holds for every section $\s' \in \Wk{s}(E_1)$ satisfying $|\jH{k}\s'| \le
c$.
\end{proposition}

\begin{proof}
The first part of the proposition is an immediate corollary to
Propositions~\ref{prop:jet-map} and \ref{prop:palais-axiom5}. The second
part of the proposition follows from estimate in
Proposition~\ref{prop:palais-axiom5}. 
\end{proof}

\begin{remark}
  \label{rem:linear-operators} 
The restriction $s \geq
  2n+4+k$ in Proposition~\ref{diffop-est} can be relaxed to $s\geq k$ when
  $D$ is a linear operator, and the estimate assumes the form
$ \| D\s \|_{s-k} \leq C \|\s \|_{s}$.
When $D$ is nonlinear, the term $\|D\s\|_{s-k}$ involves products of $\s$ and
its derivatives, and estimating these expressions uses
Lemma~\ref{lemma:loc-poly-estimate}, which assumes  $s -k \geq 2n+4$.
\end{remark}

\begin{examples}
\label{example:operator-order}
As examples, we consider the contact order of some basic operators that we
need later.
Let $M$ be a contact manifold with contact form $\n$ and characteristic
vector field $T$.
\begin{enumerate}
\item[(i)] Lie differentiation with respect to the Reeb vector field $T$
$\displaystyle
      \Lie_{T}: \Wk{s}(M) \to \Wk{s-2}(M) \,:\, f \mapsto T(f) 
$
is a differential operator of contact order~$2$.
To see this,  we work locally, using an adapted coordinate chart $\phi: U \to
  \R^{2n+1}$ as in Section~\ref{notation}.
Note that $\phi_* T = T_0 = \pd{}{x^{2n+1}}$. But since $T_0= [X_1,X_{n+1}]$, the
Lie bracket of the two horizontal vector fields, it follows that $T$ has
contact order~2.

\item[(ii)] The exterior derivative operator
$\displaystyle
  d:  \Wk{s} (\Lambda^{p}M) \to \Wk{s-2}(\Lambda^{p+1}M)
\st \a \mapsto d\a
$
 is a differential operator of contact order~$2$.
We again work locally. First consider the case $p=0$, where $\a = f$, for
$f$ a scalar function on $M$. In this case,
\[
	    d f = \sum_{k=1}^{2n} X_k(f) dx^k + T_0(f) \n_0 \,.
\]
Since the term $T_0(f)\n_0$ has contact order~2, and all other terms depend
only on the horizontal 1-jet of $f$, we see that $d\a$ has contact order~2.
For $p>0$, $\a$ can be expressed in the form
\[
      \a =        \sum_I   f_I\,dx^I  + \sum_K g_K\, \n_0 \wedge dx^K  \,,
\]
where
$
       dx^I =   dx^{i_1}\wedge \dots \wedge dx^{i_p}
$ and 
$         \n_0\wedge dx^K = \n_0 \wedge dx^{k_1} \wedge \dots \wedge dx^{k_{p-1}}$,
where summation ranges over indices $I$ and $K$ of the forms
 $i_1 < i_2 < \dots <i_p<2n+1$ and  $k_1 < k_2<\dots <k_{p-1} < 2n+1$.
Hence,
\[
      d\a = \sum_I d f_I \wedge dx^I + \sum_K d g_K\wedge\n_0 \wedge dx^K 
		       + \sum_K g_K d\n_0\wedge dx^K  \,,
\]
which is clearly of contact order~2.

\item[(iii)] Exterior differentiation followed by projection onto
$\Lambda^{*}H^*$ is a differential operator of contact order~1:
\[
        \dH:  \Wk{s}(\Lambda^p M) \rightarrow \Wk{s-1}(\Lambda^{p+1}H^{*})
\,:\, \a \mapsto  d\a - \eta \wedge (T\inter d\a)\,.
\]
To see this, let $\a$ be given in local coordinates as in (ii). Then an
easy computation gives
\[
\dH\a 
      = \sum_I \dH f_I \wedge dx^I 
		       + \sum_K g_K d\n_0\wedge dx^K  \,.
\]
But $\dH f_I =  \sum_{k=1}^{2n} X_k(f_I) dx^k$, which is of contact
order~$1$. 

\item[(iv)] Now suppose $s \geq 2n+4$, let $E\to M$ be a smooth fibre bundle, and
let $\beta$ be a fixed, smooth $p$-form on the total space $E$. The
restriction on $s$ ensures that the Folland-Stein space $\Wk{s}(E)$ is
well-defined. (Since $\Lambda^{p}M\to M$ is a vector bundle, the
restriction $s-2 \geq 2n+4$ is not required to define
$\Wk{s-2}(\Lambda^{p}M)$.)  We claim that the differential operator
\[
	 P_{\beta}: \Wk{s}(E)
\to \Wk{s-2} (\Lambda^{p}M) \,:\,\s \mapsto \s^*\beta
\]
is a
differential operator of contact order~$2$.
To see this, we work in adapted product coordinates $(x,y)\in
\R^{2n+1}\times\R^m$ on $E$ (see Section~\ref{notation}). A section $\s$ of
$E$ is then represented locally by a function $y = f_{\s}(x)$, which (by
the Sobolev lemma) is at least of class $C^{1}$. The form $\beta$ can be
written
\begin{equation}
\label{eq-beta-local}
    \beta = \sum_{|I|+|K| = p} A_{I,K}(x,y)\, dx^I \wedge dy^K
 +\sum_{|I|+|K| = p-1} B_{I,K}(x,y)\, \n_0 \wedge dx^I \wedge dy^K
 \,
\end{equation}
where $I$ and $K$ are the obvious increasing multi-indices. Notice that the
pullback $\s^*\beta$ is obtained by setting $y = f_{\s}(x)$ and expanding.  But
for a function $f_{\s}$, the differential $f_{\s}\mapsto df_{\s}$ is a differential
operator of contact order~2, since it involves derivatives in the $T_0$
direction.
 Applying this observation to each
term in the above expansion of $\beta$ shows that  $\s \mapsto \s^*\beta$ has
contact order~2.

\item[(v)] For $s \geq 2n+4$, $E$ and $\beta$ as above,
\[
	 P_{H,\beta}:\Wk{s}(E) \to \Wk{s-1} (\Lambda^{p}H^*)
	  \,:\, \s \mapsto    \piH \s^*\b 
\]
is a differential operator of contact order~$1$.
 Using the expansion~\eqref{eq-beta-local} above, we have the local identity
\[
    \n\wedge \s^{*}\beta = 
\sum_{|I|+|K| = p} A_{I,K}(x,y)\, \n_0\wedge dx^I \wedge df_{\s}^K \,.
\]
Now observe that
\[
 df_{\s} = 
  \sum_{j=1}^{n} X_j(f_{\s}) dx^j + X_{n+j}(f_{\s}) dx^{n+j} 
+ T_{0}(f_{\s})\n_0\,.
\]
It follows that the local expression for  $\n\wedge \s^{*}\b$ only
involves first derivatives of $f_{\s}$ with respect to
$X_1,\dots,X_{2n}$. Hence, $\s \mapsto  \n\wedge\s^{*}\b$  has contact
order~1. To see that $\s \mapsto \piH\s^*\b$  has contact order $1$,
notice that by Equation~\eqref{eq:piH-formula}, $\piH\s^*\b =
T\inter (\n\wedge\s^*\b)$.
Since
 $\a\mapsto T\inter\a$ is a smooth  map of vector bundles, it 
 preserves $\Wk{s}$-spaces. Hence,  the composition
$\s \mapsto \n\wedge\s^*\b \mapsto T \inter \left( \n\wedge \s^*\b \right)$ also has
contact order~1.
\end{enumerate}
\end{examples}

\subsection{Rumin's complex}
\label{sec:rumin}

In \cite{Rumin:1994}, Rumin constructed a novel resolution 
\[
0 \hookrightarrow \R \longrightarrow \cRQ^0
\stackrel{\dQ}{\longrightarrow} \cRQ^1 \stackrel{\dQ}{\longrightarrow}
\dots
\stackrel{\dQ}{\longrightarrow} \cRQ^n
\stackrel{\DQ}{\longrightarrow} \cRQ^{n+1}
 \stackrel{\dQ}{\longrightarrow} \cRQ^{n+2} \stackrel{\dQ}{\longrightarrow} \cdots
\stackrel{\dQ}{\rightarrow} \cRQ^{2n+1} \longrightarrow 0 \,
\]
of the constants on a contact manifold. In this section, we give a brief
sketch Rumin's construction.

For $n<p \leq 2n+1$,  $\RQ^{p}$ denotes the subbundle of $\Lambda^{p}M$
given by
\[
		\RQ^{p} := \{ \b\in \Lambda^{p} M \st \n \wedge \b = 0
		\text{ and } d\n \wedge \b = 0 \} 
\]
and for $0 \leq p \leq n$, $\RQ^{p}$ denotes the quotient bundle
 $\RQ^{p} =
\Lambda^{p}(M)/I^{p}$, where $I^{0}=0$ and
\[
  I^{p} : = \{\n \wedge \a + d\n \wedge \b \st \a \in \Lambda^{p-1} M, \b \in \Lambda^{p-2}M
			   \} \text{ for $0<p \leq n$.}
\]
Also note that for $0 \leq p\leq n$, $\RQ^p$ can be written as a quotient bundle of $\Lambda^p H^*$:
\begin{equation}
\label{eq:tau-quotient}
          \tau:\Lambda^p H^*\to \RQ^p =\Lambda^p H^* / (d\n\wedge
\Lambda^{p-2} H^*) \,.
\end{equation}
Let $\cRQ^{p} = C^{\infty}(\RQ^{p})$.
Then  $\cRQ^{0} = \W^0(M)$,  and since we make the identification $H^* = \RQ^{1}$ with
the annihilator of $T$, 
\[
     \cRQ^1 \ \equiv \{\a\in \W^1(M) \st T\inter\a = 0 \}\,.
\]

The linear differential operators $\dQ$ and $\DQ$ are induced by the
exterior derivative operator on forms.  Let $\b \in \cRQ^{p}$ be any section
of $\RQ^{p}$ There are three cases to consider:
\begin{enumerate}
\item[(i)] 
For $p>n$,  set $\dQ\b = d\b$. It is easy to see that $d\b$ is a section of
$\RQ^{p+1}$.
\item[(ii)] 
For $p<n$, set $\dQ\b = \piQ(d\widetilde{\b})$, where $\widetilde{\b} \in \W^{p}(M)$ is any $p$-form with
$\piQ\widetilde{\b} = \b$ and $\piQ: \Lambda^{p} M \to \RQ^{p}$ denotes the
quotient map. It is not difficult to check that $\dQ\b$ is independent of
the choice of $\widetilde{\b}$.
\item[(iii)] For $p=n$, set $\DQ(\b) = d\widetilde{\b}$, where
  $\widetilde{\b}\in \W^n(M)$ is an $n$-form satisfying the conditions
  $\piQ\widetilde{\b} = \b$ and $d\widetilde{\b} \in \cRQ^{n+1}$. Rumin shows that a
  form $\widetilde{\b}$ satisfying these conditions exists and that
  $d\widetilde{\b}$ is independent of the choice of $\widetilde{\b}$.
\end{enumerate}
Rumin also shows that $\dQ$ and $\DQ$ are linear differential operators, with
$\dQ$ of contact order~1 and  $\DQ$ of contact order~2. Using the
star operator, Rumin proves that $ \cRQ^{k} $ is dual to 
$\cRQ^{2n+1-k}$ and that the adjoint operators satisfy the identities
\[
   \deltaQ = (-1)^k *\dQ*
	  \text{ for $ k \neq (n+1)$ and } \DQs = (-1)^{n+1} *\DQ* \,.
\]
Thus $\deltaQ$ has contact order~1 and $\DQs$ has contact
order~2. 
The next proposition then follows from Proposition~\ref{diffop-est} above.
\begin{proposition}[Rumin]
Let $(M^{2n+1}, \n)$ be a compact contact manifold with
adapted metric $g$ and the associated
complex $(\cRQ^*, \dQ)$. Then the following estimates hold
\begin{align*}
&\begin{cases}
  \|\dQ \a\|_{s-1} \le c_s \|\a\|_s\,, & \text{for $k \neq n$}\\
  \|\DQ \a\|_{s-2} \le c_s \|\a\|_s\,, & \text{for $k = n$;}
\end{cases}
\quad
\text{ and }
\quad
\begin{cases}
\|\deltaQ \a\|_{s-1} \le c_s \|\a\|_s \,, & \text{for $k \neq n+1$}\\
\|\DQs \a\|_{s-2}    \le c_s \|\a\|_s \,, & \text{for $k=n+1$.}
\end{cases}
\end{align*}
\end{proposition}

Rumin\cite{Rumin:1994} establishes a Hodge theory for this complex.
The Laplace operators $\DeltaQ:\cRQ^k\to\cRQ^k$  are defined as follows

\begin{equation}
\label{def:laplace}
\DeltaQ =
\begin{cases}
 (n-k) \dQ \deltaQ + (n-k+1)\deltaQ \dQ \,, &\text{for }0 \le k \le (n-1) ,
\\
( \dQ \deltaQ)^2 + \DQs \DQ \,, & \text{for $k=n$},\\
\DQ \DQs + (\deltaQ \dQ)^2  \,,  &\text{for $k=n+1$},\\
 (n-k + 1) \dQ \deltaQ + (n-k)\deltaQ \dQ  \,, &\text{for $(n+2) \le k \le (2n + 1)$}.
\end{cases}
\end{equation}

\begin{theorem}[Rumin]
\label{thm:rumin}
Let $M$ be a contact manifold with adapted metric $g$. The
Laplace operators $\DeltaQ$  are maximal hypoelliptic,
and the following estimates are satisfied for $\a \in \cRQ^k$:
\[
\begin{cases}
     \|\a\|_{s+2} \le c_s \|\DeltaQ \a\|_s + \|\a\|_0
	       & \text{for  $k \ne n, (n+1)$,} \\
     \|\a\|_{s+4} \le c_s \|\DeltaQ \a\|_s + \|\a\|_0
	       & \text{for $k = n, (n+1)$.}
\end{cases}
\]
\end{theorem}
\begin{remark}
Rumin only establishes  the estimates stated above for the case $s=0$ (see
\cite[page 290]{Rumin:1994}).
However, standard regularity theory yields the estimate for general 
$s > 0$. For a self contained proof, one may refer to \cite{BD2}.
\end{remark}

In the following corollary, parts (i) and (ii) were explicitly stated in
\cite{Rumin:1994}; the commutation relations follow from the definitions;
the other parts follow from the hypoelliptic estimates by standard
arguments.

\begin{corollary}[Rumin]
\label{cor:rumin}
Let $(M^{2n+1}, \n)$ be a compact contact manifold with
adapted metric $g$ and the associated
complex $(\cRQ^*, \dQ)$; let $\DeltaQ$
be the associated Laplacian.

\noindent (i) The cohomology of the complex is  finite dimensional and
represented by $\DeltaQ$-harmonic forms.

\noindent(ii)
There exist operators $\GQ, \HQ : \cRQ^k\to \cRQ^k$ such that 
\[
 \text{Id} = \GQ \DeltaQ + \HQ = \DeltaQ \GQ + \HQ\,,
\]
inducing the   orthogonal decompositions:
\[
 \cRQ^k =
       \ker\DeltaQ \oplus \image \DeltaQ  
      =  \ker  \DeltaQ \oplus  
      \left( \image  \dQ \oplus \image \deltaQ \right).
\]
In particular, each  $\a \in \cRQ^k$ has a \emph{Hodge decomposition}
\[
     \a = \begin{cases}
  \HQ(\a) \oplus  (n-k)\GQ\dQ\deltaQ (\a) \oplus (n-k+1) \GQ\deltaQ\dQ(\a)
&\text{for $k <n$}\\
  \HQ(\a) \oplus \GQ(\dQ\deltaQ)^2(\a) \oplus \GQ\DQs\DQ(\a)
&\text{for $k=n$.}
	  \end{cases}
\]
\noindent(iii) The following commutation relations are satisfied: \\
For $P$ any of the operators $\dQ$, $\deltaQ$, $\DQ$, or $\DQs$,  
\[
              P \HQ = \HQ P = 0 \,;
\]
for $\a \in \cRQ^k$,
\[
 \dQ\GQ (\a) = \begin{cases}
\displaystyle\left( \frac{n-k-1}{n-k+1} \right) \GQ\dQ (\a)&\text{for } k \leq n-2\\
\displaystyle\frac{1}{2} \GQ (\dQ\deltaQ\dQ)(\a) & \text{for }k=n-1\,,
\end{cases}
\]
and 
\[
 \deltaQ\GQ (\a) = \left( \frac{n-k+2}{n-k} \right) \GQ\deltaQ (\a)
\,,\text{ for } k \leq n-1 \,.
\]
\noindent(iv)
For $\a \in \cRQ^k$,
$
\begin{cases}
 \| \GQ \a\|_{s+2} \le c_s \| \a \|_s  &\text{for $k \ne n, (n+1)$}\\
 \| \GQ \a\|_{s+4} \le c_s \| \a \|_s  & \text{for $k = n, (n+1)$}
\end{cases}
$
and  $\|\HQ (\a)\|_s \le c_s \| \a \|_s $.

Moreover, since the space of harmonic forms is finite
dimensional, $\|\a\|_s \le c_s \|\a\|_0$ for all $\a \in
\ker \DeltaQ$. In particular, $\HQ(\a)$ is of class $C^{\infty}$ for all
$\a \in \Wk{s}(\RQ^k)$.
\end{corollary}


\subsection{Characterization of contact vector fields}
\label{sec:cont-vec-fields}
In this section, we present a characterization of the closed subspace
$\Wk{s}_{cont}(TM)\subset \Wk{s}(TM)$ of $\Wk{s}$ contact vector fields in
terms of the Hodge decomposition of the Rumin complex.  We use this
characterization in Section~\ref{sec:param} to give a parameterization of
the space of contact diffeomorphisms near the identity by contact vector
fields near $0$.

We begin by recalling a few well known facts about contact vector fields.
Recall that a smooth vector field $X$ is called a \emph{contact vector field} if and only if
$\Lie_{X}(\n) = 0$ mod $\n$ (or, equivalently, $\piH(\Lie_{X}(\n))=0$).
Write $X$ in the form $X = X^0 T + \XH $  where $\XH \in H $ and $T$ is
the Reeb vector field, and use the identity
\[
	   \Lie_X (\n) = X \inter d \n + d (X \inter \n)
		       = \XH  \inter  d \n + d X^0 
\]
to conclude  that $X$ is a contact vector field if and only if it satisfies
the standard identity 
\[
     - \XH \inter d\n + T(X^0)\n = d X^0 \,.
\]
Recalling that  $\RQ^1 = H^*$ and  $\dQ f = \piH(df)$, for $f \in
C^{\infty}(M) = \cRQ^0$, we can  express this characterization in terms of
the Rumin complex as follows: $X$ is a contact vector field if and only if
it satisfies the identity
\begin{equation}\label{contact.vf}
    	 \dQ X^0 = - \XH  \inter d\n  \, .
\end{equation}
Next observe that $d\n$ defines a vector bundle
map
\[
     (\;)^{\flat}:  TM  \to H^{*} \st X \mapsto X^{\flat} = X \inter d\n
\]
whose restriction to $H\subset TM$ is an
isomorphism between the contact distribution and its dual space, and let
\[
    (\;)^{\sharp}: H^{*} \to H \st \phi \mapsto \phi^{\sharp}
\]
denote its inverse. The map
\[
    C^{\infty}(M,\R) \to C_{cont}^{\infty} \,:\, g 
    \mapsto X_g = g T - (\dQ g)^{\sharp}\, 
\]
is an isomorphism between the space of smooth functions on $M$ and the
space of smooth contact vector fields.  This map then extends to an
isomorphism between the weighted spaces:
\begin{equation}
\label{eq:gen-function}
    \Wk{s+1}(M) \to \Wkcont{s} \,:\, g \mapsto X_g = g T - (\dQ 
    g)^{\sharp} \, .
\end{equation}
The only new result here is the gain of one derivative in the
  Folland-Stein spaces, which follows easily from the two inclusions $g \in
  \Wk{s}$ and  $\dQ g \in  \Wk{s}$.

We can now express the condition for $X$ to be a contact vector field in
terms of the harmonic decomposition. Notice, in particular, the additional regularity in $X^0$,
the Reeb component of $X$. (The restriction $s \geq 2n+4$ below ensures
that $X$ is  of class $C^1$.)

\begin{lemma}
\label{lemma:char}
A vector field $X \in \Wk{s}(TM)$, $s \geq 2n+4$, is contact if and only if
it satisfies each of the following three conditions
\begin{subequations}
\label{eq:cont}
\begin{align}
\label{eq:cont-i}
\tag{a}    &X^0 =   \HQ(X^0) -  (n+1)\, \GQ\deltaQ(X\inter d\n)\\
\label{eq:cont-ii}
\tag{b}    & \HQ(X\inter d\n) = 0\\
\label{eq:cont-iii}
\tag{c}    &\begin{cases}
	 \dQ(X\inter d\n) = 0 &\text{for $n>1$}\\
	 \DQ(X\inter d\n) = 0 &\text{for $n=1$\,.}
    \end{cases}
\end{align}
\end{subequations}
Moreover, if $X \in \Wk{s}(TM)$ is a contact vector field, 
then $X^0 \in \Wk{s+1}(M)$.
\end{lemma}

\begin{proof}
Suppose that $X \in \Wk{s}(TM)$. 
Applying the Hodge decomposition to Equation~\eqref{contact.vf} shows 
that $X$ is a contact vector field if and only if it satisfies each
of the three conditions 
\[
     \HQ(\dQ X^0 + X\inter d\n) =0, 
\quad
    \deltaQ(\dQ X^0 + X\inter d\n) = 0,
\text{ and }
\begin{cases}
    \dQ(\dQ X^0 + X\inter d\n) =0 &\text{if $n>1$}\\
    \DQ(\dQ X^0 + X\inter d\n) =0 &\text{if $n=1$}
    \end{cases}
\]
The middle equation  is equivalent to
\[
       \GQ \deltaQ(\dQ X^0 + X\inter d\n)
=\frac{1}{(n+1)}( \GQ\DeltaQ X^0) + \GQ\deltaQ(X\inter d\n) = 0\,,
\]
from which the conditions~\eqref{lemma:char} follow.
Finally, suppose that $X\in \Wk{s}(TM)$ is a contact vector field.
By \eqref{eq:cont-i} and Corollary~\ref{cor:rumin},
the function $X^0$ is an element of $\Wk{s+1}(M)$.
\end{proof}


 \section{The topological group of contact  diffeomorphisms}
\label{sec:top-group-str}

Let $\diff{s}\subset \Wk{s}(M,M)$, $s \geq 2n+4$ denote the subspace of
$\Wk{s}$-diffeomorphisms of $M$.  It is well known that the space of
$C^1$-diffeomorphisms is an open subset of the space $C^1(M,M)$ of
$C^1$-maps. Moreover, since $s\geq 2n+4$, there is a continuous inclusion
$\Wk{s}(M,M)\subset C^1(M,M)$. It follows that $\diff{s}$ is an open subset of
$\Wk{s}(M,M)$.  It is, therefore, a smooth infinite dimensional manifold;  but it is
not a group  because $\diff{s}$ is not closed under composition
(see Remark~\ref{rem:invariance-1}).

Let $\cdiff{\infty} \subset \diff{s}$ denote the subspace of $C^{\infty}$
contact diffeomorphisms. By definition, the space of $\Wk{s}$ contact
diffeomorphisms of $M$ is the closed subspace $\cdiff{s} :=
\overline{\cdiff{\infty}} \subset\diff{s}$.  We show in this section that
$\cdiff{s}$ is closed under composition and inversion, and that both
operations are continuous. Consequently, $\cdiff{s}$ is a topological group
for $s \geq 2n+4$.  In Section~\ref{sec:param}, we prove that $\cdiff{s}$
is a smooth Hilbert manifold (see Theorem~\ref{thm:contact-smooth}).  Our
approach in this section and in the following section parallels the
treatment of the full diffeomorphism group given by Ebin~\cite{Ebin:1970}.

\subsection{Continuity of composition}

To prove that composition is continuous, it is sufficient to work locally.
Consider open domains $D\Subset \R^{2n+1}$ and $ \widetilde{D} \Subset
\R^{2n+1}$. By the Sobolev lemma (see
Remark~\ref{rem:sobolev-fibre-bundle}), there is a continuous inclusion
$\Wk{s}(D,\R^{2n+1}) \subset C^{1}(D,\R^{2n+1})$. Consequently, the
topological subspace $\cdiffD{s}$ of $C^1$ contact diffeomorphisms $f$ with
$f(\bar{D}) \subset \widetilde{D}$ is well defined.

\begin{proposition}
\label{prop:comp-estimate}
Let $s \ge 2n+4$, and let $f  \in \cdiffD{s}$ and
$g \in \Wk{k}(\widetilde{D},\R^m)$ for $k \le s$.
Let $D' \Subset D$ be an open set.
Then the restriction to $D'$ of the composition $g \circ f$ is an
element of $\Wk{k}(D',\R^m)$.  Moreover, the map
\[
     \mu: \cdiffD{s} \times \Wk{k}(\widetilde{D},\R^m) \to \Wk{k}(D',\R^m)
\,:\, (f,g) \mapsto g \circ f
\]
is continuous.
\end{proposition}

\begin{proof}
Our proof mimics the proof of Ebin \cite[Lemma~3.1]{Ebin:1970} in the case
of diffeomorphisms of a manifold.  It proceeds by induction on $k$.

For $k=0$, we first note that $\|g\circ f\|_{D',0} < \infty$ for any $D'
\Subset D$. Since $f$ is a $C^1$ diffeomorphism on $D$, its Jacobian
determinant $Jf$ is continuous and bounded below on $\bar{D'}$ by a
positive constant; and (by the change of variables formula for integration)
\[
   \|g\circ f\|_{D',0} =  \int_{D'} (g\circ f)^2 dV_{0}
	= \int_{f(D')} g^2\, (1/Jf\circ f^{-1}) \, d V_0 <
    \infty \,.
\]
To prove continuity at $(f,g)$, choose $\epsilon>0$.  
We will show that $\| g'\circ f' - g \circ f \|_{D',0}^2 < 4\epsilon$ for
$(f',g')$ sufficiently near $(f,g)$.  To see this, choose $\delta>0$ such that
 $\max_{\bar{D'}}(\ 1/J(f'), 1/J(f)) < \epsilon/\delta$
whenever  $\|f' - f\|_{D,s} < \delta$. Also choose a
smooth function $g^{\infty}$ on the closure of $\widetilde{D}$ such that
\mbox{$\| g - g^{\infty} \|_{\widetilde{D},0}^2 < \delta$}, and set
\[
	     M = \left(
	    \max_{x\in \widetilde{D}}
	     \sum_{i=1}^{2n+1} \left| \pd{ g^{\infty}(x) }{x^i} \right|^2
      \right) \,.
\]
 Then
\begin{align*}
\|g \circ f - g'\circ f' \|^2_{D',0}
\leq& \,
	 \|  g         \circ f - g^{\infty} \circ f \|^2_{D',0}
    +    \| g^{\infty} \circ f - g^{\infty} \circ f'\|^2_{D',0} \\
     & + \| g^{\infty} \circ f' - g         \circ f'\|^2_{D',0}
      +  \| g          \circ f' - g'        \circ f'\|^2_{D',0} \\
\leq& \, \frac{\epsilon}{\delta} \int_{\widetilde{D}} | g - g^{\infty} |^2 dV_0 +
     M \cdot  \|f - f'\|_{D,0}^2 \\
    & + \frac{\epsilon}{\delta} \int_{\widetilde{D} }
	      | g - g^{\infty} |^2 dV_0
      + \frac{\epsilon}{\delta} \int_{\widetilde{D} } |g - g'|^2 dV_0\\
\leq& \,  2 \epsilon + M \cdot  \|f-f'\|^2_{D,0} +
       \left( \frac{\epsilon}{\delta} \right) \|g -
       g'\|^2_{\widetilde{D},0} \,.
\end{align*}
Now let
$\displaystyle
     \delta' = \min\left( \delta,
       \frac{\epsilon}{M } \right)
$.
Then the last line is bounded by $4 \epsilon$ provided that $f'$ and $g'$ satisfy
the inequalities
\[
     \|f - f'\|^2_{D,0} < \delta' \text{ and }
     \|g - g' \|^2_{\widetilde{D},0} < \delta' \,.
\]

Assume that for some $k \geq 0$ the proposition holds for all $D$, $D' \Subset
D$, and $\widetilde{D}$.
We first show that $g\circ f$ is an element of $\Wk{k+1}(D')$, $k+1 \leq
s$, for all $g \in \Wk{k+1}(\widetilde{D},\R^m)$.  To do this, we need only
show that $X_{j} (g\circ f)$ is in $\Wk{k}(D')$ for $1 \leq j \leq 2n$,
where $X_j$ are the horizontal vector fields defined in
Section~\ref{notation}.  Begin by observing that since $f$ is a contact
diffeomorphism, it's derivative $f_*$ respects the contact distribution on
$\R^{2n+1}$:
\begin{equation}
\label{eq:defA}
          f_{H,*}: H_{x}\to H_{f(x)} \,:\, X_j(x) \mapsto f_*(X_j(x)) =
          \sum_{i=1}^{2n+1} A_j^i(x) X_i(f(x)),\quad 1 \leq j \leq 2n \,,
\end{equation}
where $A_j^i \in \Wk{s-1}(D,\R)$ depend continuously on $f$.
This permits us to compute as follows using the chain rule:
\begin{equation}
\label{eq:basis-change}
           X_{j}(g\circ f) = dg( f_*(X_j)) = \sum_{i=1}^{2n} A_j^i \cdot
           X_i(g)\circ f\,.
\end{equation}
By the induction hypothesis,  $X_i(g) \circ f \in \Wk{k}(D'')$ for
any open set $D''$ such that $D' \Subset D'' \Subset D$.
Since $s-1 \geq 2n+3$, we can apply Lemma~\ref{lemma:algebra} to the
products $A_j^i \cdot (X_i(g)\circ f)$
conclude that $A_j^i \cdot (X_i(g)\circ f) \in \Wk{k}(D')$, which in turn
shows that
 $X_{j}(g\circ f) \in \Wk{k}(D')$.
To complete the induction step, we have to prove continuity of
composition.  First note that if $g'$
is near $g$ in $\Wk{k+1}(\widetilde{D})$ then $X_i(g')$ is near $X_i(g)$
in $\Wk{k}(\widetilde{D})$.  Now choose a fixed open set $D''$ with $D'
\Subset D'' \Subset D$.  By the induction hypothesis, if $f'$ is near $f$
in $\Wk{s}(D)$, then $X_i(g') \circ f'$ is near $X_i(g)\circ f$ in
$\Wk{k}(D'')$.  But then by Lemma~\ref{lemma:algebra} it follows that
$\displaystyle X_{j}(g'\circ f') = \sum_{i=1}^{2n} (A_j^i)^{'} \cdot
\left( X_i(g') \circ f' \right) $ is near $\displaystyle X_{j}(g\circ f) =
\sum_{i=1}^{2n} A_j^i \cdot \left( X_i(g) \circ f \right) $ in
$\Wk{k}(D')$.
\end{proof}

\begin{corollary}
\label{cor:comp1}
Let $M$ be a compact contact manifold of dimension $2n+1$, and let $N$ be a
smooth manifold of dimension $m$.
Then
 the composition map
\[
  \mu: \cdiff{s} \times \Wk{k}(M,N) \to \Wk{k}(M,N) \,:\, (F,G) \mapsto G \circ F
\]
is continuous for $2n+4 \leq k \leq s$. In case $N = \R^m$, the map is
continuous for $2n+4 \leq s$ and $0 \leq k \leq s$.
\end{corollary}

\begin{proof}
  In case $N=\R^m$, choose $s \geq 2n+4$ and $0 \leq k\leq s$. Continuity
  of $\mu$ follows easily from the previous proposition.  We next consider
  the case where $N$ is an arbitrary smooth manifold and the definition
  of $\Wk{k}(M,N)$ requires that $2n+4 \leq k$. Fix $F \in \cdiff{s}$ and
  $G\in \Wk{k}(M,N)$. First note that the restriction $2n+4 \leq k \leq s$
  ensures that the spaces $\cdiff{s}$ and $\Wk{k}(M,N)$ are both well
  defined and that both $F$ and $G$ are of class $C^1$.  Choose adapted
  atlases $\{(\phi_{\a},U_{\a},D_{\a})\}$ and
  $\{(\widetilde{\phi}_{\a},\widetilde{U}_{\a},\widetilde{D}_{\a})\}$ for
  $M$, and charts $\{\psi_{\a},V_{\a}, B_{\a}\}$ for $N$, such that for all
  $\a$,
\[
	     F(\overline{U_{\a}}) \subset \widetilde{U}_{\a} \,,
      F_{\a}(\overline{D_{\a}})  \subset \widetilde{D}_{\a} \,,
\text{ and } G(\overline{ \widetilde{U}_{\a}})    \subset V_{\a}\,,
\]
where $F_{\a} = \widetilde{\phi}_{\a} \circ F \circ \phi_{\a}^{-1} \in
\Wk{s}(U_{\a},\widetilde{U}_{\a})$ and $G_{\a} = \psi_{\a} \circ G \circ
\widetilde{\phi}_{\a}^{-1} \in \Wk{k}(\widetilde{U}_{\a},\R^m) $.  Set
$H_{\a} = G_{\a}\circ F_{\a} (=\psi_{\a} \circ G\circ F \circ
\phi^{-1}_{\a})$.  By Proposition~\ref{prop:comp-estimate}, $H_{\a} \in
\Wk{k}(\widetilde{U}_{\a},\R^m)$ for all $\a$, showing that $G\circ F$ is
an element of $\Wk{k}(M,N)$.  To prove continuity of $\mu$, consider the
open neighbourhoods of $F$, $G$, and $H=G\circ F = \mu(F,G)$
\begin{align*}
   O(F,\epsilon,\{U_{\a}\}) &= \{ F' \in \cdiff{s}
	 \st F' (\overline{U_{\a}}) \subset  \widetilde{U}_{\a}
       \text{ for all $\a$, } \max_{\a} \| F'_{\a} - F_{\a}\|_{U_{\s},s} < \epsilon \}\,, \\
   O(G,\epsilon, \{\widetilde{U}_{\a}\} )
     &= \{ G' \in \Wk{k}(M,N)
	 \st G'(\overline{\widetilde{U}_{\a}}) \subset V_{\a}
       \text{ for all $\a$,}
	    \max_{\a} \| G'_{\a} - G_{\a}\|_{\widetilde{U}_{\s},k} <
            \epsilon \}\,,\\
   O(H,\epsilon, \{U_{\a}\} )
     &= \{ H' \in \Wk{k}(M,N)
	 \st G'(\overline{U_{\a}}) \subset V_{\a}
       \text{ for all $\a$,}
	    \max_{\a} \| H'_{\a} - H_{\a}\|_{U_{\a},k} <
            \epsilon \}\,.
\end{align*}
By definition of the topology of $\Wk{k}(M,N)$, every open neighbourhood of
$H$ contains a set of the form $O(H,\epsilon, \{U_{\a}\} )$ for
sufficiently small $\epsilon$. Moreover, by
Proposition~\ref{prop:comp-estimate}, for every $\epsilon>0$ there exists a
$\delta>0$ such that
\[
     H' = G'\circ F' \in O(H,\epsilon, \{U_{\a}\})
\]
for all $F' \in O(F,\delta, \{U_{\a}\})$ and $G' \in O(G,\delta,
\{\widetilde{U}_{\a}\})$.  Therefore,
\[
  \mu\left( O(F,\delta,\{D_{\a}\}) \times O(G,\delta,\{\widetilde{D}_{\a}\}) \right)
 \subset O(G\circ F,\epsilon, \{U_{\a}\}) \,.
\]
This completes the  proof of continuity of $\mu$.
\end{proof}

\subsection{Continuity of inversion}

The proof of continuity of inversion relies on the next lemma.
\begin{lemma}
\label{lemma:inverses}
Let $f\in \cdiffD{s}$, $s \geq 2n +4$ be a contact diffeomorphism,
with $f(D)\Subset \widetilde{D}$, and let $\widetilde{D}'\Subset
f(D)$.   Then $f^{-1} \in
\cdiff{s}(\widetilde{D}',D)$.  Moreover, the map $\iota: f' \mapsto
{f'}^{-1}$, $f' \in \{f'\in \cD_{cont}^{s}(D,\widetilde{D}),\st
\widetilde{D}'\subset f'(D)\}$, is continuous at $f'=f$.
\end{lemma}

\begin{proof} 

  Since $s \geq 2n+4$, $f$ is a $C^1$ contact diffeomorphism on every
  compact subset of $D$.  Hence, $f^{-1}$ is a $C^1$ contact diffeomorphism
  on every open, compactly contained, subset of $f(D)$.  Let $A:D \to
  GL(2n)$ be the matrix valued function defined by $ A_j^i$, where $A_j^i
  \in \Wk{s-1}(D)$ are defined as in the proof of
  Proposition~\ref{prop:comp-estimate}.  Because the process of inverting
  $A$ only involves multiplication, addition, and division of functions in
  $\Wk{s-1}$, and because $s-1> 2n+3$, we can invoke
  Lemmas~\ref{lemma:algebra} and \ref{lemma:division} to conclude that
  $A^{-1}\in \Wk{s-1}( D', GL(2n) )$ for all $D'\Subset D$.

Next observe that Equation~\eqref{eq:basis-change} with $g$ replaced by $g
\circ f^{-1}$ assumes the form
\[
           X_{j}(g) = \sum_{i=1}^{2n} A_j^i \cdot
           X_i(g\circ f^{-1})\circ f\,.
\]
Multiplying by $B = A^{-1}$ and composing with $f^{-1}$ then yields the formula
\begin{equation}
\label{eq:recursion}
     X_{k}(g\circ f^{-1})= \sum_j (B^{j}_{k} \cdot X_{j}(g) ) \circ f^{-1}.
\end{equation}
To show that $f^{-1} \in \Wk{s}(\widetilde{D}',\R^{2n+1})$ for every open
set $\widetilde{D}' \Subset f(D)$, it suffices to show that $X_{I}(f^{-1})
\in \Wk{s-k}(\widetilde{D}')$ for every multi-index $I$ with $|I| = k$ (see
Section~\ref{notation}).  Following the argument on \cite[page
17]{Ebin:1970}, we proceed by induction on $k$ to show that for all $I$
with $|I| = k$
\begin{equation}
\label{eq:inductionI}
        X_{I}(f^{-1}) =g_{I} \circ f^{-1 } \text{ with $g_{I} \in \Wk{s-k}$.}
\end{equation}
To see that \eqref{eq:inductionI} holds for $k=1$, let $g=id_{\R^{2n+1}}$ in
\eqref{eq:recursion} to get
\[
     X_{i}(f^{-1})= \left( \sum_j (B^{j}_{i}\cdot
       X_{j}(id_{\R^{2n+1}})\right) \circ f^{-1} := g_{i}\circ f^{-1}\,,
\]
and recall that $B\in \Wk{s-1}$ to conclude that $g_{i} \in \Wk{s-1}$.
Now assume that \eqref{eq:inductionI} holds for $s>k>0$, and let $I = (i,J)$, for
$J$ a multi-index with $|J|=k$. Then applying \eqref{eq:recursion} gives
\begin{equation}
\label{eq:defgI}
       X_{I}(f^{-1}) = X_i ( g_{J}\circ f^{-1} )
=   \left( \sum_j (B^{j}_{i} \cdot X_{j}(g_J)) \right) \circ f^{-1} := g_{I}
\circ f^{-1}.
\end{equation}
Since $B^j_i \in \Wk{s-1}$, $s-1\geq 2n+3$ and $X_{j}(g_J) \in \Wk{s-k-1}$,
we can invoke Lemma~\ref{lemma:algebra} to conclude that $X_{I}(f^{-1})$ is
of the form $g_{I} \circ f^{-1}$, for $g_{I} \in \Wk{s-k-1}$. This
completes the induction step.  We now know that $X_{I}(f^{-1}) = g_{I}
\circ f^{-1}$ with $g_{I}\in \Wk{s-k} \subset \Wk{0}$, for all $I$ with
$|I| \leq s$. But since $f^{-1}$ is of class $C^1$, the composition
$X_{I}(f^{-1})=g_{I}\circ f^{-1}$ is also in $\Wk{0}$ for all $I$ with $|I|
\leq s$. Hence, $f^{-1}$ is in $\Wk{s}(\widetilde{D}',D)$ for all
$\widetilde{D}'\Subset f(D)$.

To prove
continuity of the map $\iota: f' \mapsto {f'}^{-1}$,
we first show by finite induction that the map
\[
f' \mapsto g'_{I}:= X_{I}({f'}^{-1}) \circ f' \in  \Wk{0}(\widetilde{D}',D)
\]
depends continuously on $f'\in \cD_{cont}^{s}$ for all
$I$ with $|I|\leq s$.
Let $k=1$, and note that by definition of $A$ (see
Equation~\eqref{eq:defA}), the assignment $f' \mapsto A' \mapsto B' =
{A'}^{-1} \in \Wk{s-1}$ depends continuously on $f'\in
\cD_{cont}^{s}$. Hence, $g'_i\in \Wk{s-1}$ depends continuously on $f'\in
\cD_{cont}^{s}$.  Now assume that $f'\mapsto g'_{J} \in \Wk{s-k}$ depends
continuously on $f'$ for all $J$, $|J| = k$. Set $I= (i,J)$. Then
$X_i(g'_{J})\in \Wk{s-k-1}$ depends continuously on $f'$. Hence by
\eqref{eq:defgI}, $g'_{I}\in\Wk{s-k-1}$ depends continuously on $f'$,
completing the induction step.

Thus, for any multi-index $I$ with $|I|\leq s$,
\begin{align*}
      \|X_{I}({f'}^{-1}-f^{-1})\|_{\widetilde{D}',0} &=
      \|g'_{I}\circ{f'}^{-1}-g_{I}\circ f^{-1})\|_{\widetilde{D}',0}\\
 &\leq 
\|g'_{I}\circ {f'}^{-1} - g_{I}\circ {f'}^{-1}\|_{\widetilde{D}',0} +
\|g_{I}\circ {f'}^{-1} - g_{I}\circ f^{-1}\|_{\widetilde{D}',0} \,.
\end{align*}
Because the map $f'\mapsto {f'}^{-1}$ is continuous in the $C^1$-topology,
by making making $\|f' - f\|_{\widetilde{D}',s}$ sufficiently small, we
ensure that $\|X_{I}({f'}^{-1}-f^{-1})\|_{\widetilde{D}',0}$ is 
arbitrarily small for all $I$ with $|I|\leq s$. This concludes the proof
of continuity of $\iota$.
\end{proof}

\begin{theorem}
\label{thm:group-top}
Let $s\geq 2n+4$.
Then $\cdiff{s}$ is a topological group with group multiplication
\[
  \mu: \cdiff{s} \times \cdiff{s} \to \cdiff{s} \,:\, (F,G) \mapsto G \circ F
\]
and group inverse
\[
  \iota: \cdiff{s} \to \cdiff{s} \,:\, F \mapsto F^{-1} \,.
\]
\end{theorem}

\begin{proof}
Continuity of $\mu$ is contained in Corollary~\ref{cor:comp1}.  
Continuity of $\iota$ follows from Lemma~\ref{lemma:inverses} by an
argument similar to the one used in the proof of
Corollary~\ref{cor:comp1}. In brief, for fixed $F$, choose adapted atlases
$\{(\phi_{\a}, U_{\a},D_{\a})\}$, and $\{(\phi'_{\a},U'_{\a},D'_{\a})\}$
such that $ D'_{\a} \subset F_{\a}(D_{\a})$ for all $\a$. Then by
Lemma~\ref{lemma:inverses}, for all $\epsilon>0$ there there is a
$\delta>0$ such that for all $G \in O(F,\delta)$ we have $G^{-1} \in
O(F^{-1},\epsilon)$.
\end{proof}


\section{The smooth manifold of contact diffeomorphisms}
\label{sec:param}

In this section, we obtain a local coordinate chart for the set of contact
diffeomorphisms in a neighbourhood of the identity. As a corollary, we 
show that for $s \geq 2n+4$, the topological manifold $\cdiff{s}$ is a
smooth submanifold of the smooth manifold
$\diff{s}$ of $\Wk{s}$ diffeomorphisms of $M$.

\subsection{The smooth manifold of contact diffeomorphisms}
We begin by constructing a smooth atlas for $\diff{s}$. Our construction is
based on the following well known parameterization of smooth
diffeomorphisms near the identity diffeomorphism by smooth vector fields.
Fix a $C^{\infty}$ metric
adapted to the contact structure (see Section~\ref{notation}), and let
$\exp: TM \rightarrow M$ denote its exponential map. Recall that $\exp$ is the
$C^{\infty}$ map defined by the formula
\[
               \exp(X) := \gamma_{x,X}(1)
\]
where $\gamma_{x,X}:\R\to M$ is the unique geodesic curve with
$\gamma_{x,X}(0) = x$, $\gamma'(0) = X\in T_{x}M$. Also recall that for
$|X|$ sufficiently small, $|X|$ is equal to the Riemannian distance between
$x$ and $\exp(X)$.  Next consider the map $\chi$ from the space of
$C^1$-vector fields to the space of $C^{1}$-maps
\[
        \chi: C^{1}(TM) \to C^{1}(M,M) \st X \mapsto \FX\,,
\]
where $\FX$ is the $C^1$ map defined by composition
\begin{equation}
\label{eq:defFX}
            \FX\,:\,  M \stackrel{X}{\longrightarrow} TM
             \stackrel{\exp}{\longrightarrow} M \,.
\end{equation}
By compactness of $M$, there is a number $r>0$ such that any two points at
distance less then $r$ apart are joined by a unique length minimizing
geodesic. Let $B_rM \subset TM$ denote the bundle over $M$ of tangent
vectors of length less than $r$. It is a well known theorem in
Riemannian geometry that $\chi$ restricts to a diffeomorphism between the
space $C^{1}(B_r M)$ of $C^1$-vector fields of length less than $r$ and the
open set

\[
         \{ F \in  C^{1}(M,M) \st \dist_M(x,F(x)) < r \text{ for all $x\in
           M$} \} \,,
\]
where $\dist_M$ denotes Riemannian distance.  
\begin{proposition}
\label{prop:exp-1}
For $s \geq 2n+4$, the map $\chi$ restricts to a smooth map
\[     \chi^{s}: \Wk{s}(TM) \to \Wk{s}(M,M) \,:\,
		     X \mapsto \FX = \exp \circ X 
\]
on the space of $\Wk{s}$-vector fields.  Moreover, there is an open
neighbourhood $\cU \subset \Wk{2n+4}(TM)$ of the zero section, such that
for all $s \geq 2n+4$, $\chi^s$ restricts to a diffeomorphism
\[
    \chi^{s}: \cU^{s} \equiv \cU \cap \Wk{s}(TM)  \to \diff{s}
\]
between $\cU^{s}$ and a neighbourhood of the identity $id_{M}\in \diff{s}$.
\end{proposition}

\begin{proof}
First observe that a  map $F:M\to M$ can be viewed as a section of the
trivial fibre bundle
$
             \pi_{M}: M \times M \to M \st (x,y) \mapsto x 
$
and that $\exp$ defines a smooth map of fibre bundles
\[
   \widetilde{\exp}: TM \to M \times M \,:\, X \mapsto (\pi(X),\exp(X))
\]
where $\pi:TM\to M$ is projection onto the base point.  Smoothness of
$\chi^{s}$ then follows
from Proposition~\ref{prop:palais-axiom5}.  
Let   $U_rM = \{ (x,y) \in M\times M \st \dist_{M}(x,y) <r \}$. 
Then $\widetilde{\exp}: B_r M \to U_r M$ is a smooth
fibre bundle isomorphism. By Proposition~\ref{prop:palais}, the restriction
of $\chi$ to $\Wk{s}$-spaces
\[
          \chi^{s}: \Wk{s}(B_r M) \to \Wk{s}(U_rM) 
\]
is therefore a diffeomorphism.  To complete the proof, let $\cU$ denote the
preimage of the open set $ \Wk{2n+4}(U_r M)\cap \diff{2n+4}$ under
$\chi^{2n+4}$.\end{proof}

\begin{remark}
\label{rem:diff-atlas}
We can use $\chi^s$ to construct a smooth atlas for smooth manifold
$\diff{s}$.  Let $G$ be a smooth diffeomorphism of $M$. By
Corollary~\ref{cor:comp1}, that composition on the left with $G$ gives a
smooth diffeomorphism of $\diff{s}$. Consequently, the map
\[
    \chi^{s}_{G} :=   L_{G}\circ \chi^{s}: \cU^s \to \diff{s} \,:\, X \mapsto G\circ F_X
\]
is a local diffeomorphism. Since the set of $C^{\infty}$ diffeomorphisms of
$M$ is dense in $\diff{s}$, letting $G$ range over all diffeomorphisms
gives a smooth atlas. Since composition of smooth maps is smooth,
smoothness of the transition functions is automatic.
\end{remark}

\begin{remark}
\label{rem:TDiff}
We recall the standard construction of the tangent bundle
$\pi_{\diff{}}:T\diff{s}\to \diff{s}$ (see
\cite{Ham,Palais:1968} for background). To get a tangent vector to
$\diff{s}$ at  $F_0 \in \diff{s}$, let $t\mapsto F_t$ be a
smooth curve in $\diff{s}$ passing through $F_0$. Then $F_t$ is a smooth
family of $C^1$ diffeomorphisms, so we can
differentiate pointwise with respect to $t$ to obtain the vector field
\[
       X: M \to TM \st x \mapsto  \dot{F}_0(x) \in TM_{F_0(x)} 
\]
over $F_0$, where we have used the notation $\dot{F}_0 (x):=\left. \frac{d
    F_t(x)}{d t}\right|_{t=0}$. Conversely, given a vector field $X\in
\Wk{s}(M,TM)$ with $\pi \circ X = F_0$, observe that for small $t$ the
composition
\[
          F_t : M \stackrel{ t X}{\longrightarrow} TM 
                  \stackrel{\exp}{\longrightarrow} M
\]
is a smooth family in $\diff{s}$ and that $X = \dot{F}_0$.
An easy way to obtain the manifold structure on the total space $T\diff{s}$ is
to note that by Corollary~\ref{cor:left-composition}, composition with the
projection map $\pi:TM \to M$ induces a smooth map
$L_{\pi}^{s}:\Wk{s}(M,TM) \to \Wk{s}(M,M)$. Let $T\diff{s}$ be the preimage
of $\diff{s}$ under $L_{\pi}^{s}$ and let $\pi_{\diff{}} = L_{\pi}^{s}$.
\end{remark}

\subsection{Characterization of contact diffeomorphisms}

Recall that a diffeomorphism $F$ is a contact diffeomorphism if and only if
the pullback $F^*\n$ is a multiple of $\n$. Since this condition is
equivalent to the equation $\piH (F^*\n) = 0$, where $\piH: T^*M\to H^*$ is
the quotient map, the space of $\Wk{s}$-contact diffeomorphisms near the
identity is parameterized by the subspace
\begin{equation}
\label{def:cV}
     \cV^s \equiv \{ X \in \cU^s \st  \piH(\FX^*\n) = 0 \} \subset \Wk{s}(TM)\,.
\end{equation}
It is convenient to view the equation $\piH(\FX^*\n)=0$ in terms of the
Rumin complex. Notice that  $\RQ^1 \equiv H^*$, hence $\FX\in \cU^s$ is a
contact diffeomorphism if and only if it satisfies the equation
\[
                 \piQ(\FX^*\n) = 0 \,. 
\]
This suggests studying  the non-linear differential operator
\[
            X \mapsto \piQ(\FX^*\n) 
\]
in more detail.  Our  goal is to show 
that
\begin{equation}
\label{eqn:quadratic-term}
            \piQ(\FX^*\n) = \piQ\Lie_{X}\n + \piQ\circ\Quad_{\n}(X)\,,
\end{equation}
where $\Lie_{X}$ denotes Lie differentiation with respect to the vector
field $X$ and $\Quad_{\n}(X)$ is a smooth differential operator that
vanishes to second order as $X\to 0$. Part (ii) of the next proposition
shows that $X \mapsto \piQ(\FX^*\n)$ is a smooth differential operator of
contact over $1$; that it has the form of
Equation~\eqref{eqn:quadratic-term} is a corollary to
Lemma~\ref{lem:exp-q-form}.

\begin{proposition}
\label{prop:exp2}
The following maps are smooth (non-linear) differential operators for all $s
\geq 2n+4$:
\begin{enumerate}
\item[(i)] $\displaystyle \Wk{s}(TM) \to \Wk{s-2}(T^{*}M)
       \,:\, X  \mapsto \FX^*\n$,

\item[(ii)] $\displaystyle  \Wk{s}(TM) \to  \Wk{s-1}(R^{1})
       \,:\, X  \mapsto \piQ\FX^*\n$,

\item[(iii)] $\displaystyle \Wk{s}(TM) \to \Wk{s-1}(\RQ^2)
       \,:\, X  \mapsto \dQ \left( \piQ\FX^*\n \right) = \piQ (d \FX^*\n) = \piQ \FX^*(d\n)$, 
          for $n>1$.
\end{enumerate}
\end{proposition}

\begin{proof}
  View $\FX:M\to M$ as a section of the trivial bundle $M\times M \to
  M$. Since $s \geq 2n+4$, Example~\ref{example:operator-order} parts (iv)
  and (v) apply to yield (i) and (ii).  To prove (iii), apply
  Example~\ref{example:operator-order} (v) to the smooth 2-form $d\n$ to
  conclude that the map
\[
           X  \mapsto   \piH (\FX^*d\n)
\]
is a smooth differential operator of contact order 1. Next observe that
$\piQ = \tau \circ \piH$, where $\tau$ is the quotient map given by
\eqref{eq:tau-quotient}.  Finally, since $\dQ\b = \piQ(d\b)$, the
composition
\[
     X \mapsto \piH(\FX^*(d\n)) \mapsto 
    \tau( \piH(\FX^*d\n)) =  \piQ \FX^*(d\n) = \piQ d(\FX^*\n)
\]
is also a smooth differential operator of contact order 1.
\end{proof}

To show that $\piH(\FX^*\n)$ has the form of
Equation~\eqref{eqn:quadratic-term}, we work locally, choosing an adapted
atlas $\phi_{\a}: U_{\a}\to \R^{2n+1}$ for $M$ and a collection of open
sets $W_{\a}\Subset U_{\a}$ covering $M$ as in Section~\ref{notation}.  By
compactness of $M$, there is a constant $c>0$ such that $\exp(x,X) \in
U_{\a}$ for all $x \in \overline{W_{\a}}$, all $X\in TM_x$, with $|X|<c$,
and all $\a$.

Let $X$ be a $C^1$ vector field with $|X|<c$. Fix a chart, say $\phi_{\a}$, and
set $U = U_{\a}$ and $W =W_{\a}$. To simplify notation, we adopt the
Einstein summation conventions, letting Roman indices range from $1$ to
$2n+1$. Then there exist smooth functions $B^{k}_{ij}(x,X)$ (locally
defined) on $TM$ such that
\begin{equation}  \label{exp-map}  \exp^{k}(x,X) = x^{k} +
 X^{k} + B^{k}_{ij}(x,X)X^{i}X^{j}\, . 
\end{equation} 
This follows simply from the second order Taylor's formula with
 integral remainder for the exponential map. Indeed, for fixed
 $X \in T_{x}M$, let $\g(t) =
 \exp(x,tX)$ be a  geodesic. Then 
 \begin{equation} 
\label{eqn:exp-quad}  
\g^{k}(1) = \g^{k}(0) + \dot \g^{k}(0) - \int_{0}^{1}(1-t)
 \Gamma^{k}_{ij}(\g(t)) \dot\g^{i}(t) \dot\g^{j}(t) dt\, ,
 \end{equation} 
 where $\Gamma^{k}_{ij}$ are the Christoffel symbols, and we have used the
 geodesic equation $\ddot{\g}^{k} + \Gamma^{k}_{ij}\dot{\g}^{i}\dot{\g}^{j}=0$.
 Let $y = \exp(x, X)$.  Since $\g(t) = \exp(x,t X) = y(x, tX)$, then
\[
\dot\g^{i}(t) = \dot y^{i}(x, tX) = \pd{y^{i}(x, tX)}{X^{j}}  X^{j} \,,
\]
and  this becomes
\begin{equation*} 
\g^{k}(1) = \g^{k}(0) + \dot \g^{k}(0) - \int_{0}^{1}(1-t)
\Gamma^{k}_{ab}(\exp(x,t X)) \pd{y^{a}}{X^{i}}(x,tX)
\pd{y^{b}}{X^{j}}(x,tX) X^{i} X^{j} dt\, ,
\end{equation*} 
whence
\begin{equation} 
B^{k}_{ij}(x,X) = - \int_{0}^{1}(1-t)
 \Gamma^{k}_{ab}(\exp(x,tX)) \pd{y^{a}}{X^{i}}(x,tX)  \pd{y^{b}}{X^{j}}(x,tX) dt\, .
 \end{equation}

\begin{lemma}
\label{lem:exp-q-form}
Let $\psi$ be a smooth $q$-form on $M$ and
choose a coordinate patch $U = U^{\a}$, with $W = W_{\a} \Subset U$. Let
$c>0$ be chosen so that $\exp(x,X) \in U$ for all $x \in \overline{W}$ and
all $X \in T_xM$ with $|X|<c$.
Then  there are (locally defined) smooth fibre bundle maps
\[
    \Quad^1_{ij}: \left.   BM\right|_{W} \to \left.\Lambda^{q}M \right|_{W}
\text{ and }
    \Quad^2_{ij}: \left. BM\right|_{W} \to \left.\Lambda^{q-1}M \right|_{W} \,,
\]
where  $BM  = \{ X \in TM \st |X| < c\}$,
such that  for any $C^1$ vector field $X:M \to BM \subset TM$ the equation
\[
	\FX^{*} \psi = \psi + \Lie_{X} \psi
       +  \Quad^1_{ij}( X )\,  X^{i}X^{j}
       +  \Quad^2_{ij}( X ) \, \wedge X^{i} d X^{j}
\]
is satisfied  on all of $W$.
\end{lemma}

\begin{proof}
Begin with the special case of a $0$-form $u \in C^{\infty}(M,\R)$.
Then $\FX^*u (x)  = u \circ \exp(x,X)$, and applying Taylor's formula with
integral remainder to the function $f(t) = u(x + t X + t^2 B_{ij}(x,tX)X^i
X^j)$ and setting $t=1$ yields the formula
\[
   (u \circ \exp)(x,X) = u(x) + \Lie_X u(x) + \Quad_{ij}(x,X)X^{i}X^{j}\,,
\]
for $\Quad_{ij}(x,X)$ smooth functions  on ${\left.BM\right|}_{W}$, such that
$\Quad_{ij} = \Quad_{ji}$.
Next consider the special case  $\psi = dx^{k}$, and compute as follows, 
using what we have just proved:
\begin{align*}
   \FX^*(d x^{k}) &=
    d \left( \exp^*_{X} x^{k} \right)\\
    &= d\left( x^{k} + \Lie_{X}( x^{k} )
	+ \Quad_{ij}(x,X) X^{i}X^{j}\right)\\
    &= d x^{k} + \Lie_{X} dx^{k} + d\left( \Quad_{ij}(x,X)
    X^{i}X^{j} \right)\\
    &= d x^{k} + \Lie_{X} dx^{k}
      +  d\left(\Quad_{ij}(x,X)\right) X^{i}X^{j}
      +  2 \Quad_{ij}(x,X) X^{i}d X^{j} \\
    &= d x^{k} + \Lie_{X} dx^{k} +  \Quad^1_{ij} X^{i}X^{j}
       +  \Quad^{2}_{ij}(x,X) X^{i}d X^{j}
\end{align*}
for  $\Quad^{1}_{ij} =
\pd{\Quad_{ij}}{x^{k}}dx^{k} $, 
$\Quad^{2}_{ij} = \pd{\Quad_{ik}}{X^{j}} X^{k} + 2 \Quad_{ij} 
= \pd{(\Quad_{ik}X^{k})}{X^{j}}  +  \Quad_{ij} $.
Because every $p$-form can be expressed as a linear combination of 
products of terms as above,  the general result follows easily by induction.
\end{proof}

\begin{remark}
\label{rem:quad}
Henceforth, we will use the notation
\[
    \Quad_{\psi}(X) := \FX^*(\psi) -\psi - \Lie_X\psi  
\]
to denote the non-linear part of the pull-back $\FX^*\psi$.
 The Lemma states that in local coordinates
\[
   \Quad_{\psi}(X) =   \Quad^1_{ij}(X) \, X^{i}X^{j}
       +   \Quad^{2}_{ij}(X) \wedge \, X^{i}d X^{j} \,,
\]
where $\Quad^1_{ij}$ and $\Quad^2_{ij}$ are smooth functions on
${\left.BM\right|}_{W} \subset TM$, which depend on the smooth form $\psi$.
\end{remark}

\subsection{Parameterization of contact diffeomorphisms}
The condition for $\FX$ to be a
contact diffeomorphism is the vanishing of the one-form $\FX^* \n$ mod
$\n$.  Remark~\ref{rem:quad} applied to $\psi=\n$ and the identity
$\piQ\n=0$ show that $\FX$ is a contact diffeomorphism if and only if
\[
                 \piQ \Lie_{X}\n + \piQ \Quad_{\n}(X) = 0\,.
\]
Since $\Quad_{\n}(X)$ vanishes to second order at $X=0$, the linearization
of this equation is
\[
                \piQ \Lie_{X}\n =0\,,
\]
i.e., the condition that $X$ be a contact vector field.
This suggests using the implicit function theorem in Banach spaces to
construct a parameterization of the space of contact diffeomorphisms near
the identity by the space contact vector fields near zero.

We are going to construct a smooth map between
Hilbert spaces of the form
\[
    \Phi:\Wk{s}(TM) \to \Wk{s+1}(M) \times \cH^{s}
    \,:\, X \mapsto g_X \oplus \gamma_X \,,
\]
where $\cH^{s}$ is a second Hilbert space (to be determined), such
that
\begin{enumerate}
\item[(i)] $\FX$ is a contact diffeomorphism if and only if $\gamma_X =0$,
\item[(ii)] the derivative of  $\Phi$ is invertible at the origin.
\end{enumerate}
By the
inverse function theorem, $\Phi$ is locally
invertible and the map
\[
	 g  \mapsto \chi^{s}\left(\Phi^{-1}(g,0) \right)
\]
gives a smooth parameterization of the contact diffeomorphisms in
$\cdiff{s}$ near the identity by real valued functions in $\Wk{s+1}(M)$
near zero.

A natural guess for the map $\Phi$ is
\begin{equation}
\label{eq:guess}
	   X = X^0\,T + \XH  \mapsto (X^0 , \piQ \FX^{*}\n )\,,
\end{equation}
for, as we have already observed, $\FX$ is contact if and only if $\piQ
\FX^{*}\n = 0$, and $X^0$ parameterizes contact vector fields (see
Section~\ref{sec:cont-vec-fields}).  Unfortunately, this map is not
invertible.  Indeed, its linearization at the origin involves a
differential operator that loses too many derivatives.  

The trick to circumventing this difficulty is to exploit some  hidden
smoothness in the Hodge decomposition of one-forms in the Rumin complex.
For smooth data, choose $\Phi$ to be of the form
\[
   \Phi: \begin{cases}
     C^{\infty}(TM) &\rightarrow
     C^{\infty}(M,\R) \oplus \range(\deltaQ) \oplus \ker (\deltaQ)
		     \subset
     C^{\infty}(M,\R) \oplus C^{\infty}(M,\R) \oplus \cRQ^{1} \\
     \quad X &\mapsto   g_{X}  \oplus \a_{X} \oplus  \omega_{X} \,.
	 \end{cases}
\]
When $n>1$,  the Hodge theory shows that the projection
$\piQ\a\in\cRQ^1$ of a general one-form $\a \in \W^1$
has the decomposition
\begin{align*}
     \piQ\a 
&= \GQ \left\{ (n-1)\,\dQ \deltaQ  + n\, \deltaQ \dQ \right\}\piQ \a + \HQ\piQ \a\\
&= (n-1)\,\GQ\dQ \deltaQ\piQ\a + n\, \GQ\deltaQ \piQ d\a + \HQ\piQ\a.
\end{align*}
Applying the commutation relations of Corollary~\ref{cor:rumin}(iii), gives
\[
   \piQ\a  =
  (n+1)\, \dQ \GQ \deltaQ\piQ \a + n\, \GQ \deltaQ \piQ d\a + \HQ\piQ \a \, .
\]
This, together with the identity $d\FX^*\n = \FX^*d\n$ give the following
 decomposition of $\piQ\FX^* \n$:
\begin{equation}
\label{eq:FXeta}
     \piQ\FX^* \n =
(n+1)\, \dQ \GQ \deltaQ (\piQ\FX^* \n) 
 + n \GQ \deltaQ (\piQ\FX^* d\n) + \HQ (\piQ\FX^* \n)\,.
\end{equation}

Referring now to our natural guess \eqref{eq:guess}, we reassemble it using
Equation~\eqref{lemma:char}(a) for $X^0$ and the identity \eqref{eq:FXeta}
to define the map $\Phi$ by the formulas
\begin{equation}
\begin{cases}
  g_{X}     &=\,  -(n+1)\, \GQ \deltaQ (\XH  \inter d\n) + \HQ (X^0)\,, \\
 \a_{X}     &=\; (n+1)\,\GQ \deltaQ (\piQ \FX^* \n) \,,\\
\omega_{X}  &=\;  n\, \GQ \deltaQ \dQ  (\piQ \FX^* \n) +  \HQ (\piQ \FX^*\n )  \,.
\end{cases}
\end{equation}
Indeed, we observe that $\piQ\FX^*\n = \dQ\a_X + \w_{X}$ by
Equation~\eqref{eq:FXeta}, and in the case where $X$ is a contact vector
field, then $g_X = X_0$ by Lemma~\ref{lemma:char}. In the case $n=1$,
we have to adjust the map $\Phi$ to reflect the Hodge decomposition at
$\cRQ^n = \cRQ^1$:
\[
   \piQ\FX^* \n = \GQ (\dQ \deltaQ)^2 (\piQ\FX^* \n) +
   \GQ \DQs \DQ (\piQ\FX^* \n) + \HQ (\piQ \FX^* \n) \, .
\]
Applying the commutation relation \ref{cor:rumin}(iii) to the $0$-form
$\deltaQ(\piQ\FX^*\n)$, yields the identity
\[
   \piQ\FX^* \n = 2\dQ \GQ\deltaQ(\piQ\FX^* \n) +
   \GQ \DQs \DQ (\piQ\FX^* \n) + \HQ (\piQ \FX^* \n) \, .
\]
Because $n+1=2$, the formulas for $g$, $\a$, and $\omega$ become
\begin{equation}
\label{eq:phi-n=1}
\begin{cases}
 g_{X}     &=\,  - (n+1)\GQ\deltaQ (\XH  \inter d\n) + \HQ  (X^0)   \,,\\
 \a_{X}    &=\;  (n+1)\, \GQ\deltaQ (\piQ\FX^* \n)  \,,\\
\omega_{X} &=\;   \, \GQ \DQs \DQ (\piQ\FX^* \n) +  \HQ  (\piQ\FX^* \n)   \, ,
\end{cases}
\end{equation}
and only the formula for $\omega_{X}$ has changed.

\begin{remark} 
Observe that by construction $\piQ\FX^*\n = \dQ\a_{X} + \omega_{X}$. Since  $\dQ$
is injective on $\range(\deltaQ)$, $\a_{X}=0$ if and only if $\dQ\a_{X} =0$.
Consequently, $F_X$ is a contact diffeomorphism if and  only if
$\gamma_{X}= \a_{X}\oplus\omega_{X}=0\oplus 0$.
\end{remark}

\begin{remark}
In the case $n=1$, the map $\Phi$ is, roughly speaking, the same as the map
defined in \cite{B}; however, in that paper, the use of the complex
Laplacian and complex operators necessitated an additional splitting into
real and imaginary parts---roughly doubling the number of terms.
\end{remark}

\begin{proposition}
\label{prop:localcoordn}
Let $\Phi$ be the map defined above and let 
\[
   \cH^{s} := \Wk{s}(\range(\deltaQ))  \oplus \Wk{s}(\ker(\deltaQ)) 
               \subset  \Wk{s}(M) \oplus \Wk{s}(\RQ^{1} )\,.
\]
Then for  $s \geq 2 n+4$, the map $\Phi$ extends to a smooth map
\[
   \Phi: 
   \begin{cases}
     \Wk{s}(TM) &\rightarrow
		 \Wk{s+1}(M)\oplus \cH^{s}
		      
		    \\
     \quad X &\mapsto   g_{X}  \oplus \gamma_{X} := g_{X}\oplus (\a_{X} \oplus  \omega_{X})
   \end{cases}
\]
The linearization of $\Phi$ at the zero vector field is given by
\[
      d\Phi(X) = \begin{cases}
	\begin{aligned}   &\left\{  -(n+1)\, \GQ \deltaQ (\XH  \inter d \n) +
            \HQ (X^0)  \right\} \oplus
        \\
         &\left\{(n+1)\, \GQ \deltaQ (\dQ X^0 + \XH  \inter d \n) 
    \oplus \left(  n\, \GQ \deltaQ \dQ (\XH  \inter d \n) +
	     \HQ  (\XH  \inter d\n) \right)  \right\} \,, \\
                   &\qquad \text{for $n>1$}
\end{aligned} \\
\\
\begin{aligned}
&\left\{ - (n+1)\,  \GQ \deltaQ(\XH  \inter d\n) +\HQ  (X^0) \right\} \oplus
 \\
& \left\{  (n+1)\, \GQ\deltaQ (\dQ X^0 + \XH  \inter d \n) \oplus 
 (\GQ \DQs \DQ (\XH  \inter d \n) +  \HQ  (\XH  \inter
   d\n) ) \right\}\,,  \\
&\qquad \text{ for $n=1$.}
\end{aligned}
\end{cases}
\]
Moreover, the linearization of $\Phi$ is invertible with inverse given by
\[
   g \oplus  (\a \oplus  \omega) \mapsto
		  (g + \a )\, T  + (-\dQ g + \omega)^{\sharp}
\]
where $\sharp$ is the isomorphism from horizontal one-forms to
horizontal vector fields induced by the
two-form  $d\n$.
\end{proposition}
The proof relies on two lemmas.
\begin{lemma}
\label{lemma:a}
Let $F\in \Wk{s}(M,M)$, $s \geq 2n+4$, be a diffeomorphism, and let $\b \in
\W^{k}$ be a smooth $k$-form, $k \le n$. Then the form
$\piQ F^{*}\b$ lies in $\Wk{s-1}(\RQ^k)$.
\end{lemma}
\begin{proof}
This is an immediate corollary to
Example~\ref{example:operator-order} (v) and the observation that for
$q\leq n$,
$\piQ =\tau\circ\piH$ where $\tau$ is the quotient map \eqref{eq:tau-quotient}.
\end{proof}

\begin{lemma}
\label{lemma:Dregularity}
Let $n=1$ and let  $F:M \to M$ be a $\Wk{s}$ (possibly not contact)
diffeomorphism,  and let $s\geq 6 ( = 2n+4)$.
Then $\DQ \piQ F^{*}\n$ is in $\Wk{s-2}(\RQ^2)$.
\end{lemma}
\begin{proof}
  Recall the definition of the operator $\DQ$ (see
  Section~\ref{sec:rumin}).  For any $\a\in\cRQ^1 \subset \W^{1}(M)$,
  $\DQ\a = d(\widetilde{\a}+f\n)$, where $\widetilde{\a}\in \W^{1}(M)$ is
  any form such that $\piQ(\widetilde{\a})=\a$ and $f\in \cRQ^{0}$ is the
  unique function such that $\n\wedge d(\widetilde{\a} + f \n) = 0$.  Let
  $\Lambda$ be the linear isomorphism
\[
      \Lambda: \cRQ^{0} \to \cRQ^{3} = \W^3(M) \,:\, h \mapsto h \n\wedge d\n
      \,.
\]
Note that $\Lambda$ is defined by a smooth vector bundle isomorphism;
it, therefore, extends to an isomorphism between
the spaces  $\Wk{s}(M)$  and $\Wk{s}(\Lambda^3M)$ for all $s$.
Consider now the case $\a = \piQ F^{*}\n$. 
The condition defining $f$ is
\[
          \n \wedge (  F^*d\n + f d\n ) = 0 \,.
\]
By Example~\ref{example:operator-order}(v), $\n\wedge F^*d\n$ is in
$\Wk{s-1}(\RQ^3)$, forcing $f$ to be in
$\Wk{s-1}(\RQ^0)$. Compute as follows
\[
        \DQ(\piQ F^*\n) = d  F^*\n  + d (f\n) = F^*d\n + d f
        \wedge \n + f d\n \,.
\]
Since wedging with $\n$ kills all terms in $df\wedge\n$
involving differentiation in directions transverse to the contact
distribution, it follows that $d(f\n)$ is in
$ \Wk{s-2}(\Lambda^2M)$.  Thus $\DQ\a$ is in $\Wk{s-2}$, concluding the
proof of the lemma.
\end{proof}

\begin{remark}
The result of Lemma~\ref{lemma:Dregularity} is somewhat surprising.
Because $\piQ F^{*}\n$ is in $\Wk{s-1}$ and $\DQ$ is an operator of contact
order 2, one would expect $\DQ \piQ F^{*}\n$ only to lie in $\Wk{s-3}$.
\end{remark}

\begin{proof}[Proof of Proposition~\ref{prop:localcoordn}]
By Proposition~\ref{prop:exp-1}, 
the map $\chi^{s}:X \mapsto \FX$ is smooth; and for sufficiently small $X$, 
$\FX$ is a $\Wk{s}$-diffeomorphism.
Recall that $\piQ F^*\n$ is in $\Wk{s-1}(\RQ^1)$.
The linear operators $\dQ $, $\deltaQ$, $\DQ$, $\DQs$, $\HQ$, and
 $\piQ$ are all bounded as  maps as follows:
\begin{align*}
    &\dQ, \deltaQ : \Gamma^s \rightarrow \Gamma^{s-1}\,,
    \,\HQ : \Gamma^{s-1} \rightarrow \Gamma^s \,,
\text{ and  }
    \piQ: \Gamma^s \rightarrow \Gamma^s \,,\\
   &
\GQ:\Wk{s}(\RQ^{k})\to \Wk{s+2}(\RQ^{k})\,,
\text{  for  } k<n\,,\\
   &\DQ, \DQs : \Wk{s} \rightarrow \Gamma^{s-2}\,,
 \GQ:\Wk{s}(\RQ^{k})\to \Wk{s+4}(\RQ^{k}) 
\text{  for  } k=n\,. 
\end{align*}
and the individual terms have the following regularity properties:
\[
   \GQ \deltaQ(\piQ \FX^* \n) \in \Gamma^s, \quad \GQ \deltaQ\dQ(\piQ
     \FX^* \n) \in \Gamma^s, \quad \HQ (\piQ\FX^* \n) \in C^{\infty} \,,
\]
for $n>1$; and  for $n=1$:
\[
   \GQ \deltaQ (\piQ\FX^* \n) \in \Wk{s}, \quad
\GQ \DQs \DQ (\piQ \FX^* d\n) \in \Wk{s}, \quad \HQ(\piQ \FX^* \n) \in C^{\infty}
\,.
\]
Thus, $\Phi$ is smooth and maps between the spaces as indicated.

The only nonlinear terms in the map $\Phi$ arise from the
presence of $\FX^* \n$. The linearization of this term at the
zero vector field is $\Lie_X \n = d X^0 + \XH  \inter d\n$. When we
substitute this into the map $\Phi$, we obtain easily the
linearization and its inverse.
\end{proof}

Let $\cU \subset \Wk{2n+4}(TM)$ be an open neighbourhood of the zero
section such that $\FX$ is a $\Wk{s}$-diffeomorphism for all $X \in
\cU^{s}:=\cU\cap \Wk{s}(TM)$, $s\geq 2n+4$.  The next theorem shows that
the subset $\cV^s \subset \cU^s$ on which $\FX$ is a contact diffeomorphism
is a smooth submanifold which is smoothly parameterized by the space of
$\Wk{s}$-contact vector fields near zero.

\begin{theorem}
\label{thm:second-order}
For $\cU$ sufficiently small,  for all $s \geq 2n+4$ the set
\[
       \cV^s = \{ X \in \cU^s \st \Phi (X) = (g_X, 0, 0) \}
\]
is a smooth submanifold of $\cU^s := \cU \cap \Wkcont{s}$, smoothly
parameterized by the map
\[
\Psi : \Wkcont{s}\cap \cU^s \rightarrow \cV^s \,:\, 
        X \mapsto \Phi^{-1}( g_X \oplus 0 ) \,.
\]
Moreover, the map $\Psi$ is of the form
\[
	   \Psi(X) = X + B(X)(X,X) \, ,
\]
where $B:(\Wkcont{s} \cap \cU) \times \Wkcont{s} \times \Wkcont{s} 
	   \rightarrow \Wk{s}(TM)$ is smooth and bilinear in the last two factors.
\end{theorem}

\begin{proof}
  It suffices to prove the theorem for $s=2n+4$. That $\cV^s\subset \cU^s$
  is a smooth submanifold follows from Proposition~\ref{prop:localcoordn}
  and the inverse function theorem in Banach spaces.  To define $\Psi$, let
  $\pi$ denote the projection
\begin{equation}
\label{eqn:def-pi}
    \pi (g, \gamma) = (g, 0) \, ,
\end{equation}
in the notation in the statement of Proposition~\ref{prop:localcoordn}, and let
 $\Psi = 
\left(\Phi^{-1} \circ \pi \circ \Phi\right)|_{\Wkcont{s}\cap \cU}
$.   The
smoothness of $\Psi$ follows from the smoothness of $\Phi$. A simple
calculation shows that $d\Psi |_{0}(X) = X$, for all $X\in \Wkcont{s}$,
which (together with the inverse function theorem) shows that $\Psi$
parameterizes $\cV^s$.
The form of the operator $B$ is given by Taylor's
formula with integral remainder for smooth operators on Banach spaces:
\[
          B(X) = \int_0^1 (1-t) D^2\Psi(t X) \,dt\,.
\]
(See e.g. \cite[Theorem 3.5.6]{Ham}.)
\end{proof}

\begin{remark}
\label{rem:parm-gen-functions}
In view of the isomorphism $\Wk{s+1}(M) \simeq \Wkcont{s}$ given by
\eqref{eq:gen-function}, the map $g \mapsto \Psi(X_g)$ defines a smooth
parameterization of $\cV^{s}$ by $\Wk{s+1}$-functions in a neighbourhood of
$0\in \Wk{s+1}(M,\R)$.

\end{remark}

\subsection{The smooth structure on the space of contact diffeomorphisms}

The map $\Psi$ of
Theorem~\ref{prop:localcoordn} gives a parameterization of the subspace
$\cV^s\subset \cU^s$. We now show that this parameterization in turn
induces a smooth structure on the space $\cdiff{s}$ of all 
$\Wk{s}$-contact diffeomorphisms.

\begin{theorem}
\label{thm:contact-smooth}
Let $(M,\n)$ be a compact contact manifold. For $s \geq 2n+4$, the
space of $\Wk{s}$ contact diffeomorphisms is a smooth Hilbert manifold. 
\end{theorem}

\begin{proof}
We first show that the intersection of   $\cdiff{s}$ with a neighbourhood of
the identity  is a smooth submanifold of $\diff{s}$.
 To see this,
let $\chi^{s}: \cU^s \to \diff{s}$ be the diffeomorphism onto a neighbourhood of
the identity given in Proposition~\ref{prop:exp-1}.  By
Theorem~\ref{thm:second-order}, we can shrink $\cU^s$ if necessary so that
$\cV^s$ is a smooth submanifold of $\cU^s$. Now set $\cO^s_{\Id} =
\chi^s(\cU^s)$. Since $\chi^s:\cU^s\to \cO_{\Id}^s$ is a diffeomorphism and
\[
        \cdiff{s} \cap \cO^s_{\Id} = \chi^s(\cV^s)
\]
it follows that $\cdiff{s} \cap \cO^s_{\Id}$ is a smooth submanifold of
$\diff{s}$.

Next consider the open set $\cO_{F}^s =
F(\cO^s_{\Id})\subset \diff{s}$, where $F$ is an arbitrary $C^{\infty}$ contact
diffeomorphism.  Noting that $G \in \diff{s}$ is a contact diffeomorphism
if and only if $F^{-1}\circ G$ is a contact diffeomorphism shows that
the equality
\[
            \cO^s_{F} \cap \cdiff{s} = F(\cV^s)
\]
holds. Finally recall that composition on left with $F$ is a smooth
diffeomorphism of $\diff{s}$ (see
Remark~\eqref{rem:diff-atlas}) to conclude that 
$\cO^s_{F} \cap \cdiff{s}$ is a smooth submanifold of $\diff{s}$.

It remains only to show that every element of $\cdiff{s}$ is contained in
$\cO^s_{F}$ for some smooth contact diffeomorphism $F$. To see this, choose
any $G\in \cdiff{s}$ and let $F_{k}$ be a sequence of smooth contact
diffeomorphisms converging to $G$. Because composition and inversion are
continuous operations, $F^{-1}_k\circ G \to id_{M}$ as
$k\to\infty$. Therefore $ F^{-1}_{k}\circ G \in \cO^{s}_{\Id}$ for $k$
sufficiently large. Consequently $G$ is contained in $\cO^s_{F_k} =
F_{k}(\cO^{s}_{id_{M}})$ for $k$ sufficiently large.\end{proof}

\begin{remark}
\label{rem:atlas}
We can construct a smooth atlas for $\cdiff{s}$ as follows. 
By Theorem~\ref{thm:second-order}
and Remark~\ref{rem:parm-gen-functions}, there is an open neighbourhood
$\cO^{s+1}$ of $0 \in \Wk{s+1}(M,\R)$ such that
\begin{equation}
\label{eq:plot}
      \cchi{s+1}: \cO^{s+1} \to \cdiff{s} \st       g \mapsto \chi^{s}\circ \Psi^{s}(X_g) 
\end{equation}
is a homeomorphism  onto
$\cO^s_{\Id}\cap \cdiff{s}$. Its inverse is a coordinate chart centered at
the identity diffeomorphism. Composition with a smooth contact
diffeomorphism $G$ then yields the map
\[
        {\cchi{s+1}}_{,G}: \cO^{s+1} \to \cdiff{s} \st g \mapsto G \circ
        \cchi{s+1}(g) \,,
\] 
whose inverse is a coordinate chart centered at $G$.
The argument in the last paragraph of the proof of
Theorem~\ref{thm:contact-smooth} shows that the set of all such charts
forms a smooth atlas for $\cdiff{s}$.
\end{remark}

We next address the global topology of $T\cdiff{s}$.
Because $\cdiff{s}$ is a closed submanifold of $\diff{s}$, there is a
smooth inclusion $T\cdiff{s} \subset T\diff{s}$.  Using the fact that
$\diff{s}$ is an open subset of $\Wk{s}(M,M)$, and letting
$\Wk{s}(M,|TM|)$  denote the space of $\Wk{s}$-maps from $M$ into the
total space of $TM$ (i.e. forgetting the vector bundle structure on $TM$),
one sees immediately that the tangent bundle of $\diff{s}$ is the open
subset
\[
           T\diff{s} = \{ X \in \Wk{s}(M,|TM|) \st  \pi \circ X \in \diff{s}
           \}
\]
with bundle projection $T\diff{s} \to \diff{s} \st X \mapsto \pi \circ X$.
A standard computation with Lie derivatives applied to a one-parameter
family of contact diffeomorphisms then shows that a $\Wk{s}$-vector field
$X: M \to |TM|$ is in $T\cdiff{s}$ if and only if $F = \pi \circ X$ is
$\cdiff{s}$ and $X \circ F^{-1}\in \Wk{s}(TM)$ is a contact vector field.
As the next proposition
shows, $T\cdiff{s}$ is a trivial vector bundle:

\begin{proposition}
\label{prop:TDcont}
Let $X_g$ denote the contact vector field associated to the generating
function $g$. The map
\[
             \cdiff{s} \times \Wk{s}(M,\R) \longrightarrow T\cdiff{s} 
\st (g , F) \mapsto X_g \circ F \,,
 \]
is a continuous vector bundle isomorphism.
\end{proposition}

\begin{remark}
  Composition with a $\Wk{s}$-contact diffeomorphism is a continuous, but
  \underline{not} differentiable, operation. Consequently, the
  trivialization in Proposition~\ref{prop:TDcont} is not smooth. We discuss
  the smoothness of composition in Section~\ref{sec:smoothness-comp}.
\end{remark}

Proposition~\ref{prop:TDcont} is a corollary to a more general
construction.  Let $\pi:E\to M$ be a smooth vector bundle over $M$, and let
$|E|$ denote the total space of $E$, viewed as a smooth manifold,
forgetting its vector bundle structure.  Recall that $\Wk{s}(M, |E|)$
denotes the space of $\Wk{s}$-maps from $M$ into $|E|$. Because $\pi$ is
smooth, Corollary~\ref{cor:left-composition} applies to show that the map
\[
       L^{s}_{\pi}: \Wk{s}(M, |E| ) \to \Wk{s}(M, M) \st G \mapsto \pi\circ
       G\,.
\]
is smooth. Let $\WkD{s}(M,E) = (L^{s}_{\pi})^{-1}(\cdiff{s}) \subset
\Wk{s}(M, |E|)$ and
let 
\begin{equation}
\label{eqn:def-induced-bundle}
       \piD: \WkD{s}(M,E) \longrightarrow \cdiff{s}
\end{equation}
be the restriction of $L^s_{\pi}$ to $\WkD{s}(M,E)$. Notice that the vector
bundle structure on $E$ induces a vector space structure on the the fibres
of $\piD$, and we call $\piD:\WkD{s}(M,E)\to\cdiff{s}$ the (vector)
\emph{bundle of $\Wk{s}$-sections of $E$ over contact diffeomorphisms.}

\begin{lemma}
\label{lem:triviality}
Let $\pi: E\to M$ be a smooth vector bundle over $M$. Then for $s \geq 2n+4$,
the map $\displaystyle \Phi_E: \cdiff{s} \times \Wk{s}(E) \to \WkD{s}(M,E)
\st (F,\s) \mapsto \s \circ F $ is a continuous, vector bundle isomorphism
between $\WkD{s}(M,E)$ and  the trivial vector bundle.
\end{lemma}

\begin{proof}
  Consider first the special case where $E\to M$ is the trivial bundle $M
  \times \R^r\to M$. The diffeomorphism $\Wk{s}(M, |M \times \R^r|) \simeq
  \Wk{s}(M,M) \times \Wk{s}(M,\R^r)$ restricts to a diffeomorphism
\[
     \WkD{s}(M, M \times \R^r) \simeq \cdiff{s} \times \Wk{s}(M,\R^r)
\]
with respect to which $\Phi_{M\times\R^r}$ assumes the form
\[
 \Phi_{M\times \R^r}:
 \cdiff{s} \times 
\Wk{s}(M, \R^r) \to \WkD{s}(M, M \times \R^r) \st
 (F, \sigma) \mapsto (F, \sigma \circ F) \,,
\]
with inverse
\[
 \Phi^{-1} _{M\times \R^r}:\WkD{s}(M, M \times \R^r)\to \cdiff{s} \times
 \Wk{s}(M,\R^r) \st (F, \sigma) \mapsto (F, \sigma \circ F^{-1}) \,.
\]
Continuity of $\Phi^{-1}_{M\times\R^r}$ follows from continuity of composition with $F$
(see Corollary~\ref{cor:comp1}); and continuity of $\Phi^{-1}_{M\times\R^r}$
follows from continuity of inversion (see Theorem~\ref{thm:group-top}).

Now consider the general case.  By construction, $\Phi_E$ is bijective,
preserves basepoint, and is linear on each fibre.  To see that $\Phi_E$ is
continuous, note that since the map $ E\to M\times E \st e \mapsto (\pi(e),
e)$ is smooth, so is the induced map
\[
     \iota: \Wk{s}(E)\hookrightarrow \Wk{s}(M\times E)    
\]
defined by the formula $\iota(\sigma): x \mapsto (x,\sigma(x))$.  This
observation, together with Corollary~\ref{cor:left-composition} implies
continuity of $\Phi_E$.  It remains only to show that $\Phi^{-1}_E$ is
continuous.  Let $j:E \to M \times \R^r$ be a smooth vector bundle
inclusion into a trivial bundle, and let $s: M \times \R^r \to E$ be a
smooth vector bundle map with $s \circ j = id_E$.  Continuity of
$\Phi^{-1}_E$ is proved by expressing $\Phi^{-1}_E$ as the following
composition of continuous maps
\begin{multline*}
  \WkD{s}(M,E) \stackrel{L^s_{j}}{\longrightarrow}
 \WkD{s}(M, M\times \R^r) 
\simeq \cdiff{s} \times \Wk{s}(M,\R^r)\\
\stackrel{\Phi^{-1}_{M\times\R^r}}{\longrightarrow}
 \cdiff{s} \times \Wk{s}(M,\R^r)
\stackrel{ id_{\cdiff{s}}\times L^{s}_{s}}{\longrightarrow}
 \cdiff{s} \times \Wk{s}(M,E) \,,
\end{multline*}
concluding the proof of the lemma.
\end{proof}

We close this section with a formula for the derivative of $\cchi{s}$,
which we need in Section~\ref{sec:smoothness-comp}.  For
$g \in \cO^{s+1}$, we denote by $T_{F_g}\cdiff{s}$ the tangent space to
$\cdiff{s}$ at the contact diffeomorphism $F_g$,  and for $h \in
\Wk{s+1}(M,\R)$, we set   
\[
      Y_{(g,h)} =  D{\cchi{s+1}}_g (h) : M \to TM\,,
\]
where $Y_{(g,h)} \in T_{F_g}\cdiff{s}$;  and we set  $X_{(g,h)} = Y_{(g,h)} \circ F_{g}^{-1} \in \Wkcont{s}$.
\begin{lemma}
\label{lem:Dcchi2}
For $s \geq 2n+4$, the map
\[
       \cO^{s+1} \times \Wk{s+2}(M,\R) \to \Wkcont{s} \st
   (g,h) \mapsto X_{(g,h)}
\]
is continuous. Moreover, for every $g\in \cO^{s+1}$ and $\epsilon>0$, there
is a $\delta>0$ such that
\[
           \| X_{(g_1,h)} - X_{(g,h)} \|_{s} < \epsilon \|h\|_{s+1}
\]
for all $g_1\in \cO^{s+1}$ such that $\| g_1 - g\|_{s+1} < \delta$ and all
$h\in \Wk{s+2}(M,\R)$.
\end{lemma}
\begin{proof}
Continuity of the map is clear. The estimate is the restatement of the fact
that the derivative $D{\cchi{s}}_{,g}$ depends continuously on $g$.
\end{proof}

\subsection{Differentiability of composition}
\label{sec:smoothness-comp}

We showed in Section~\ref{sec:top-group-str} 
that composition
\[
  \mu: \cdiff{s} \times \Wk{k}(M,\R) \to \Wk{k}(M,\R) \,:\, 
          (F,u) \mapsto u \circ F
\]
is a continuous operation for $2n+4 \leq k \leq s$, but composition is not
$C^1$, as the following counterexample shows. Choose $u \in
\Wk{k}(M,\R)$ with $T(u) \notin \Wk{k}(M\,\R)$, where $T$ is the Reeb
vector field on $M$. Then the one-parameter family $F_t$ of contact
diffeomorphisms given by the flow of the Reeb vector field $T$ is a smooth
curve in $\cdiff{s}$; and differentiability of $\mu$ would imply that the limit
\[
         \lim_{t \to 0}  \frac{ \mu(u, F_t) - \mu(u,F_{0})}{t} = T(u)
\]
would be an element of $\Wk{k}(M,\R)$. But this contradicts our choice of $u$.
The next theorem shows that  we can recover smoothness by strengthening
the regularity assumption on $u$.

\begin{theorem}
\label{thm:comp-smoothness}
Let $M$ be a compact contact manifold of dimension $2n+1$ and let $N$ be a
smooth manifold. 
Then the map
\[
  \mu: \cdiff{s} \times \Wk{k+2}(M,N) \to \Wk{k}(M,N) 
          \,:\, (F,G) \mapsto G \circ F
\]
is continuously differentiable for  $2n+4 \leq k \leq s$. In case $N=\R^m$,
$\mu$ is continuously differentiable for $0 \leq k \leq s$, $2n+4 \leq s$.
\end{theorem}

\begin{proof}
Assume that the theorem holds in the special case where $N = \R^m$, with
$m$ arbitrary. Let $\iota: N \hookrightarrow \R^m$ be a closed embedding,
and let
  $U\subset \R^m$ be a tubular neighbourhood of $N$, with projection map
  $\pi:U\to N$. By Corollary~\ref{cor:left-composition}, the maps $\iota$
  and $\pi$ induce smooth maps
\[
     \widetilde{\iota}: \Wk{k}(M,N) \to \Wk{k}(M,U)
\text{ and }
     \widetilde{\pi}: \Wk{k}(M,U) \to \Wk{k}(M,N)
\]
for all $k \geq 2n+4$. 
Because $\Wk{k}(M,U)$ is an open subset of
$\Wk{k}(M,\R^m)$, by assumption, we know that
\[
          \mu: \cdiff{s} \times \Wk{k+2}(M,U) \to \Wk{k}(M,U): (F,u)
          \mapsto u \circ F
\]
is a $C^1$ map. It follows that the composition
\[
        \cdiff{s} \times \Wk{k+2}(M,N) \stackrel{id \times 
\widetilde{\iota}}{\longrightarrow}
 \cdiff{s} \times \Wk{k+2}(M,U) \stackrel{\mu}{\longrightarrow} \Wk{k}(M,U)
        \stackrel{\widetilde{\pi}}{\longrightarrow} 
   \Wk{k}(M,N)
\]
is  a $C^1$ map. 

It remains to prove the theorem in the case $N=\R^m$. Because
$\Wk{k}(M,\R^m)$ is the $m$-th fold product of $\Wk{k}(M,\R)$, we need only
prove it for $m=1$; and by Remark~\ref{rem:atlas}, it suffices to restrict to
an open neighbourhood of the identity in $\cdiff{s}$.
Now for
$\cO^{s+1} \subset \Wk{s+1}(M)$ a sufficiently small neighbourhood of $0$,
the map
\[
             \cchi{s+1}: \cO^{s+1} \longrightarrow \cdiff{s} \st g \mapsto
             F_g 
\]
is a smooth parameterization of a neighbourhood of the identity contact
diffeomorphism. (Here and in the following we set $F_g = F_{\Psi(X_g)}$, where $X_g$ denotes the
contact vector field with generating function $g$.)  
With this notation, the proof reduces to proving that the map
\[
     \mu: \Wk{k+2}(M,\R) \times \cO^{s+1} \longrightarrow \Wk{k}(M,\R) \st
          (u,g) \mapsto u \circ F_g 
\]
is $C^1$. The  next proposition completes the proof.
\end{proof}

\begin{proposition}
\label{prop:Dmu}
For $s \geq 2n+4$ and   $s \geq k \geq 0$ and for $\cO^{s+1} \subset \Wk{s+1}(M,\R)$ a
sufficiently small neighbourhood of $0$, the map
\[
     \mu: \Wk{k+2}(M,\R) \times \cO^{s+1} \longrightarrow \Wk{k}(M,\R) \st
          (u,g) \mapsto u \circ F_g \,.
\]
is $C^1$ with derivative at $(u,g)$ given by the formula
\[
        D\mu_{(u,g)}:  (v,h) \mapsto  v \circ F_g +    (X_{(g,h)} \inter
        du)\circ F_g\,,
\]
for $(v,h) \in \Wk{k+2}(M,\R)\times \Wk{s+1}(M,\R)$.
\end{proposition}

\begin{proof}
  Because $\mu$ is a map between Banach spaces, to show that it is $C^1$,
  we need only show that the two partial derivatives of $\mu$ with respect
  to the first and second variables
\begin{align*}
   D_1\mu:  \Wk{k+2}(M,\R) \times \cO^{s+1} &\to L\left( \Wk{k+2}(M,\R),
   \Wk{k}(M,\R) \right) \\
   D_2\mu:  \Wk{k+2}(M,\R) \times \cO^{s+1} &\to L\left( \Wk{s+1}(M,\R),
   \Wk{k}(M,\R) \right)\\
\end{align*}
exist and are continuous\footnote{We use the notation $L(\cH_1,\cH_2)$ to
  denote the Banach space of bounded linear maps between Hilbert spaces
  $\cH_1$ and $\cH_2$.}.
We shall obtain formulas for $D_1\mu$ and $D_2\mu$. The formula for
$D\mu_{(u,g)}(v,h)$ then follows immediately from the well-known identity
\[
        D\mu_{(u,g)}(v,h) = D_1\mu_{(u,g)}(v) + D_2\mu_{(u,g)}(h) \,. 
\]

To see that  $D_1\mu$ exists, notice that by
Corollary~\ref{cor:comp1}, $\mu$ is continuous. The map $\mu$ is linear in the
first variable and, therefore,   differentiable with respect to the
first variable, with derivative given by
$D_1\mu_{(u,g)}(v) =  v \circ F_g$.
Continuity of $D_1\mu$ is proved in Lemma~\ref{lem:D1mu} below.

We next claim that
\[
        D_2\mu_{(u,g)}(h) = ( X_{(g,h)}\inter du) \circ F_g \,.
\]
To prove the claim, first notice that because $s\geq 2n+4$, the functions $g$,
$h$, $u$, as well as the map $F_{g}$ are of class at least $C^1$. Moreover,
because $\cchi{s+1}$ is smooth, the family $t \mapsto F_{g+ t h}$ of
contact diffeomorphisms is a smooth family. Consequently, we can compute
pointwise at $x\in M$, employing the chain rule as follows:
\begin{multline}
D_{2}\mu_{(u,g)}(h) (x) = \lim_{t\to 0} \frac{ u (F_{g + t h}(x)) - u(
  F_{g}(x))}{t}
 = \left. \frac{d u\left(F_{g+t
        h}(x)\right)}{dt}\right|_{t=0}\\
\qquad = du \left( {D\cchi{s+1}}_{,g}(h)(x)
\right) = Y_{(g,h)}(x)\inter du_{F_g(x)} = (X_{(g,h)}\inter du) \circ
F_{g}(x)\,.
\end{multline} 
To prove that $\mu$ is differentiable with respect to the second
variable we need to verify the formula
\[
\lim_{h \to 0} 
   \frac{\| u \circ F_{g+h} - u \circ F_{g} - (X_{(g,h)}\inter du)\circ F_g\|_{k} }
          { \|h\|_{s+1}} = 0 \,
\]
and we need to prove continuity of $D_2\mu$. We do this in
Lemma~\ref{lem:D2mu}.
\end{proof}

\begin{lemma}
\label{lem:D1mu}
 $D_1\mu:  \Wk{k+2}(M,\R) \times \cO^{s+1} \to L\left( \Wk{k+2}(M,\R),
   \Wk{k}(M,\R) \right)$ is continuous for $s \geq 2n+4$ and all $k \geq 0$.
\end{lemma}

\begin{proof}
  By definition of continuity, we must show that for any $g_0 \in
  \cO^{s+1}$ and any $\epsilon>0$, there is a $\delta>0$ such that the
  condition
\[
        \| v \circ F_{g} - v \circ F_{g_0} \|_{k} < \epsilon
               \|v\|_{k+2}
\]
is satisfied for all $v \in \Wk{k}(M,\R)$ and  all $g \in \cO^{s+1}$ with
$\|g-g_0\|_{s+1} < \delta$.  We need only prove the estimate for $v$ a
smooth test function.

Set $h = g - g_0$. Because the map $g \mapsto F_g$ is smooth,
and $v$ is smooth, for fixed $x\in M$, the
function $\gamma_x: t \mapsto v( F_{g_0 + t h}(x))$
is a $C^1$~function of $t$. Consequently, we can compute as follows using
the chain rule:
\begin{align*}
 v \circ F_{g}(x) - v \circ F_{g_0}(x)
&
=\int_{0}^{1} \frac{d}{dt}   v(F_{g_0+t h}(x))   \,dt 
&
=  \int_{0}^{1}
                  (X_{(g_0 + t h,h)} \inter dv) \circ F_{g_0 + t h}(x) \,dt \,.
\end{align*}
Viewing $t \mapsto (X_{(g_0 + t h,h)} \inter dv) \circ F_{g_0 + t h}$ as a
continuous curve in the Hilbert space $\Wk{k}(M,\R)$ yields the inequality
\[
         \| v \circ F_{g}(x) - v \circ F_{g_0}(x) \|_{k} 
\leq \int_0^1 \left\| (X_{(g_0 + t h,h)} \inter dv) \circ F_{g_0 + t h}
\right\|_{k} \, dt \,.
\]
Hence, to prove the lemma it suffices to show that, for $\delta>0$
sufficiently small, the estimate
\begin{equation}
\label{eqn:D1mu1}
                   \| (X_{g,h} \inter dv ) \circ F_g \|_{k} < \epsilon
                   \|v\|_{k+2}
\end{equation}
holds for all $g\in \cO^{s+1}$ and $h\in \Wk{s+1}(M,\R)$ with
$\|g-g_0\|_{s+1} < \delta$ and $\|h\|_{s+1} < \delta$.
To obtain \eqref{eqn:D1mu1} first observe that since interior evaluation
\[
                \Wk{s}(TM) \times \Wk{k}(T^*M) \to \Wk{k}(M,\R) : (X,\beta)
                \mapsto X \inter \beta
\]
is a smooth bilinear map (see Proposition~\ref{prop:palais}),  the estimate 
\begin{equation}
\label{eqn:D1mu2}
           \| X\inter \beta \|_{k} \prec \|X\|_{s} \| \beta \|_{k}
\end{equation}
holds. Note also that by continuity of composition (see
Corollary~\ref{cor:comp1}) and linearity in $v$,
\begin{equation}
\label{eqn:D1mu3}
              \| v \circ F_{g_0} \|_{k} \prec \|v \|_{k}
\end{equation}
for all $v$. Continuity also shows that
 we can choose $\delta>0$ small so  that
\[
        \| v\circ F_{g} - v \circ F_{g_0} \|_{k}<1
\]
holds provided $\|v\|_{k} < \delta$ and $\|g - g_0\|_{s+1} <
\delta$. By linearity in $v$, setting $C_2 = 1/\delta$, we get the
estimate
\begin{equation}
\label{eqn:D1mu4}
                  \| v\circ F_{g} - v \circ F_{g_0} \|_{k}< C_2 \|v\|_{k}\,,
\end{equation}
provided $\|g-g_0\|_{s+1} < \delta$. Finally note that because
$\cchi{s+1}$ is smooth,  its derivative $D\cchi{s+1}$ is continuous in the
operator norm. We can therefore choose $\delta>0$ so small that
\begin{equation}
\label{eqn:D1mu5}
                \| X_{(g,h)} - X_{(g_0,h)} \|_{s} < \epsilon \|h\|_{s+1}
\end{equation}
for all $h$, provided $\|g - g_0\|_{s+1} < \delta$.  Choosing $\delta>0$ so
that all of the above estimates hold and choosing $g$ so that $\|g -
g_0\|_{s+1} < \delta$, we can then estimate as follows:
\begin{align*}
 \| (X_{(g,h)}\inter dv) \circ F_g \|_{k} &\leq
 \| (X_{(g,h)}\inter dv) \circ F_{g_0} \|_{k} +
 \| (X_{(g,h)}\inter dv) \circ F_g  -(X_{(g,h)}\inter dv ) \circ F_{g_0}
 \|_{k} \\
&\prec
 \| (X_{(g,h)}\inter dv) \circ F_{g_0} \|_{k} + \| X_{(g,h)}\inter dv \|_{k}
\prec
 \| X_{(g,h)}\inter dv \|_{k} + \| X_{(g,h)}\inter dv \|_{k}
 \\
&\prec
 \| X_{(g,h)}\inter dv   \|_{k}   \prec \| X_{(g,h)} \|_{s} \cdot
 \|v\|_{k+2}
\prec
 (\| (X_{(g_0,h)}\inter dv )  \|_{s}  + \epsilon \|h\|_{s+1}) \cdot  \|v\|_{k+2}
\\
&\prec \| h \|_{s+1} \cdot \|v \|_{k+2}
\end{align*}
The estimate \eqref{eqn:D1mu1} follows by decreasing $\delta$, if necessary,
and requiring $\|h\|_{s+1} < \delta$. \end{proof}

\begin{lemma}
\label{lem:D2mu}
For $2n+4 \leq s$, $k \leq s$, and  $\cO^{s+1}$ as in Lemma~\ref{lem:D1mu},
the derivative
 $D_2\mu:  \Wk{k+2}(M,\R) \times \cO^{s+1} \to L\left( \Wk{s+1}(M,\R),
   \Wk{k}(M,\R) \right)$ exists, is continuous, and given by the formula
\[
       D_2\mu_{(u,g)}(h) =  Y_{(g,h)} \inter  (du \circ F_g) \,,
\]
for $g \in \cO^{s+1}$, $u \in \Wk{k}(M,\R)$, and $h\in \Wk{s+1}(M,\R)$.
\end{lemma}

\begin{proof}
To prove that $D_2\mu$ exists, choose $(u,g_0) \in \Wk{k+2}(M,\R) \times
\cO^{s+1}$. We need to show that
\begin{equation}
\label{eqn:D2mu1}
    \lim_{h\to 0} \frac{ \| u \circ F_{g_0+h} - u \circ F_{g_0} 
- (X_{(g_0,h)} \inter du) \circ F_{g_0}
      \|_{k}} {\|h \|_{s+1}} = 0 \,.
\end{equation}
Choose $\epsilon>0$. We  need to find $\delta>0$ such that
\begin{equation}
\label{eqn:D2mu2}
             \| u \circ F_{g_0+h} - u \circ F_{g_0} 
- (X_{(g_0,h)} \inter du) \circ F_{g_0}
      \|_{k} < \epsilon \| h\|_{s+1}
\text{ for } \|h\|_{s+1}<\delta \,.
\end{equation}
To this end, compute as follows for $g\in \cO^{s+1}$ near $g_0$, setting
$h = g - g_0$:
\begin{align*}
 \| u \circ F_{g_0+h} &- u \circ F_{g_0} 
- (X_{(g_0,h)} \inter du) \circ F_{g}
      \|_{k} = \left\| \int_0^1   \left\{ \frac{d  (u \circ F_{g_0+th})}{dt} 
                -  (X_{(g_0,h)} \inter du) \circ F_{g_0} \right\}\, dt
            \right\|_{k}\\
&\leq \int_0^1 
\left\| \left\{ \frac{d  (u \circ F_{g_0+th})}{dt} 
                -  (X_{(g_0,h)} \inter du) \circ F_{g_0} \right\}
            \right\|_{k}
\,dt\\
&=\int_0^1 
\left\|  (X_{(g_0+ t h,h)}\inter du) \circ F_{g_0 + t h}  
         -  (X_{(g_0,h)} \inter du) \circ F_{g_0}             
\right\|_{k} \,dt
\end{align*}
Consequently, to prove \eqref{eqn:D2mu2}, it suffices to find $\delta>0$ so
that the inequality
\begin{equation}
\label{eqn:D2mu3}
  \left\|  (X_{(g,h)}\inter du) \circ F_{g}  
         -  (X_{(g_0,h)} \inter du) \circ F_{g_0}             
\right\|_{k} < \epsilon \|h\|_{s+1}
\end{equation}
holds for all $\|g-g_0\|_{s+1}<\delta$. But
by Proposition~\ref{prop:Dmu}, we can choose $\delta>0$ so that\\ 
$\| v\circ F_g - v \circ F_{g_0}\|_{k} < \epsilon \|v\|_{k+2}$ for all $v \in \Wk{k}(M,\R)$;
and with this choice of $\delta$, we can estimate as follows:
\begin{align*}
  &\left\| (X_{(g,h)}\inter du) \circ F_{g} 
         -  (X_{(g_0,h)} \inter du) \circ F_{g_0}  \right\|_{k} \\
&\leq 
 \left\|  
       (X_{(g,h)}\inter du) \circ F_{g} - (X_{g,h}\inter du) \circ F_{g_0} 
  \right\|_{k}
                +  
  \left\|
      (X_{(g,h)} \inter du) \circ F_{g_0} - (X_{(g_0,h)}\inter du) \circ F_{g_0} 
   \right\|_{k} \\
& \prec \epsilon \|X_{(g,h)}\inter du\|_{k} + 
         \|  \left( X_{(g,h)} - X_{(g_0,h)} \right)\inter du \|_{k} \\
& \prec \epsilon \|X_{(g,h)}\inter du\|_{k} + 
         \|  X_{(g,h)} - X_{(g_0,h)}\|_{s} \cdot \| du \|_{k} \,.\\ 
\intertext{Finally, using smoothness of $\cchi{s+1}$ as we did in Equation~\eqref{eqn:D1mu5}
above we can bound the last term as follows, for $\delta$ sufficiently small:}
&\ \prec  \epsilon \|X_{(g,h)}\inter du\|_{k} + (\epsilon \|h\|_{s+1})
 \|du\|_{k}
\prec \epsilon \|h\|_{s+1} \,.
\end{align*}
This concludes the proof of \eqref{eqn:D2mu2}.
Continuity of $D_2\mu$ follows from the following estimate:
\begin{align*}
      \| D_2\mu_{(g,u)}(h) &- D_2\mu_{(g_0,u_0)}(h)\|_{k}
= \| (X_{(g,h)}\inter du)\circ F_g - (X_{(g_0,h)}\inter d u_0)\circ F_{g_0}
\|_{k}\\
&\prec
\| (X_{(g,h)}\inter du)\circ F_g - (X_{(g_0,h)}\inter d u)\circ F_{g_0}
\|_{k} +
\| (X_{(g_0,h)}\inter du)\circ F_{g_0} - (X_{(g_0,h)}\inter d u_0)\circ F_{g_0}\|_{k}\\
&\prec
\epsilon \|du\|_{k} \|h\|_{s+1} + C  \| d(u-u_0)\| \|h\|_{s+1} \prec \epsilon \|h\|_{s+1}\,,
\end{align*}
which holds for all $(g,u)$ with $\| g - g_0\|_{s+1} < \delta$, $\|u - u_0\|_{k+1} < \delta$.
\end{proof}

\subsection{Some a priori estimates}
Theorems~\ref{thm:second-order} and \ref{thm:contact-smooth} state that the
nonlinear space of $\Wk{s}$ contact diffeomorphisms is a Hilbert manifold
modelled on the linear space of $\Wk{s}$ contact vector fields. In typical
applications, one would like to study the action of the space of contact
diffeomorphisms on some set of structures by comparing it with the
linearized action of contact vector fields. For this strategy to work, it
is necessary to show that the error incurred in the linearization is quadratically small in
an appropriate sense. This is the content of the next proposition, which gives
\apriori\ estimates for the quadratic error $\Psi(X) - X = B(X)(X,X)$. We
require this result in \cite{BD3} to obtain normal forms for CR~structures on
compact three dimensional contact manifolds.

\begin{proposition}
\label{prop:Psi-estimates}
For $X \in \Wkcont{s} \cap \cU^s$,
\begin{equation}
\tag{i}\label{sharper-estimate}
 \|\Psi(X) - X\|_{s} \prec \|X\|_{s}\|X\|_{s-1} \, .
\end{equation}
Moreover, for all  $ X_{1}, X_{2} \in   \Wkcont{s} \cap \cU^s$,
\begin{align}
\tag{ii}\label{cauchy-estimate}
\| ( \Psi(X_{2}) - X_{2})   - (\Psi(X_{1})  - X_{1}) \|_{s}
\prec  &\; \|X_{2} - X_{1}\|_{s-1}(\|X_{2}\|_{s} + \|X_{1}\|_{s})\\
 &\; + \|X_{2} - X_{1}\|_{s}(\|X_{2}\|_{s-1} + \|X_{1}\|_{s-1}) \, .\nonumber
\end{align}
\end{proposition}

Our proof relies on the next two lemmas. The first  is a corollary to
Lemma~\ref{lem:exp-q-form} and compactness of $M$. (See Remark~\ref{rem:quad}
for the definition of $\Quad$.)
\begin{lemma}
  \label{lem:remainder-estimate} For $c>0$ sufficiently small, the following
  estimates hold for $\psi$ a fixed smooth $q$ form. For all $X \in
  \Wkcont{s+2}$, $s\geq 2n+4$, such that $|X|<c$:
\begin{align}
\tag{i}\label{lem:remainder-1}
   \|\FX^*\psi \|_{s} &\prec
     \|\psi\|_s + \|\Lie_{X}\psi\|_s + \|X\|_s \, \|X\|_{s+2} \,,
\\
\tag{ii}\label{lem:remainder-2}
   \|(\FX^*\psi) \wedge \n \|_{s} &\prec
     \|\psi \wedge \n\|_s + \|\Lie_{X}\psi \wedge \n\|_s
      + \|X\|_s \, \|X\|_{s+1} \,,
\\
\intertext{and, for $q \leq n$,}
\tag{iii}
\label{lem:remainder-3}
  \|\piQ (\FX^*\psi) \|_{s}  &\prec
     \| \piQ\psi \|_s + \|\piQ (\Lie_{X}\psi ) \|_s
      + \|X\|_s \, \|X\|_{s+1} \,.
\end{align}
Moreover,  for $q \leq n$,
\begin{equation}
\tag{iv}\label{lem:remainder-4}
\begin{split}
   \|  \piQ \left(\Quad_{\psi} (X_1) - \Quad_{\psi} (X_2) \right) \|_{s} \prec&\quad
	 \|X_1 - X_2\|_s \, (\|X_1\|_{s+1} + \|X_2\|_{s+1})\\
&\quad +
	 \|X_1 - X_2\|_{s+1} \, (\|X_1\|_{s} + \|X_2\|_{s})\,,
\end{split}
\end{equation}
holds for any two  vector fields  $X_i$, $i=1,2$ such that $|X_i| < c$.
\end{lemma}
\begin{proof}
  Choose an adapted atlas and a constant $c>0$ as in the discussion above
  Lemma~\ref{lem:exp-q-form}.  The inequalities \eqref{lem:remainder-1} and
  \eqref{lem:remainder-2} follow from the definition of $\Quad_{\psi}$.  To
  prove \eqref{lem:remainder-3}, note that by Lemma~\ref{lem:exp-q-form},
  $\piQ \circ \Quad_{\psi}$ is a smooth differential operator of contact
  order $1$; the inequality then follows.  The inequality
  \eqref{lem:remainder-4} follows from the smooth dependence of
  $\Quad_{\psi}(X) $ on $X$.
\end{proof}

\begin{lemma}
\label{lem:DQ-estimate}
Choose $c>0$ as in the previous lemma. The following estimates are satisfied for
 any $\Wk{s+2}$ vector fields $X$, $X_i$, $i=1,2$ with $|X|<c$, $|X_i|<c$, $s
 \geq 2n+4$.
If $n=1$ then
\[
 \|\DQ \piQ(\FX^*\n -\n  - \Lie_X\n)\|_{s-2}
=  \|\DQ \piQ(\Quad_{\n}(X) )\|_{s-2}
\prec   \|X\|_{s}\,\|X\|_{s-1} 
\]
and
\begin{align*}
 \|\DQ \piQ(\Quad_{\n} (X_1) - \Quad_{\n}( X_2) )\|_{s-2}  
       \prec &\;   (\|X_1\|_{s} + \|X_2\|_{s})\cdot \|X_1 - X_2\|_{s-1}\\
       &\; + (\|X_1\|_{s-1} + \|X_2\|_{s-1})\cdot \|X_1 - X_2\|_{s}\,.
\end{align*}
If $n>1$ then
\[
 \|\dQ \piQ(\FX^*\n -\n  - \Lie_X\n)\|_{s-1}
=  \|\dQ \piQ(\Quad_{\n} (X) )\|_{s-1}
\prec   \|X\|_{s}\,\|X\|_{s-1} 
\]
and
\begin{align*}
 \|\dQ \piQ(\Quad_{\n} (X_1) - \Quad_{\n}(X_2) )\|_{s-1,loc}  
       \prec &\;   (\|X_1\|_{s} + \|X_2\|_{s})\cdot \|X_1 - X_2\|_{s-1}\\
       &\; + (\|X_1\|_{s-1} + \|X_2\|_{s-1})\cdot \|X_1 - X_2\|_{s}\,.
\end{align*}

\end{lemma}

\begin{proof}
Let $\phi = \FX^*\n  -\n - \Lie_X\n \in \Wk{s}(\Lambda^1(M)$.
By compactness of $M$, it suffices to obtain
local estimates on a coordinate patch $W\Subset U$ chosen
as in Lemma~\ref{lem:exp-q-form}. In the notation of
Lemma~\ref{lem:exp-q-form}, the one form $\phi$ can be written 
\begin{align*}
    \phi  
       =&\;  \Quad^1_{ij}(X) \,  X^{i}X^{j}
       + \Quad^2_{ij}(X) \,  X^{i} d X^{j} \\
       =&\;
       \Quad_{ij,k}(x,X(x)) X^{i}(x) X^{j}(x) dx^{k} +
	     \Quad_{ij}(x,X(x)) X^{i}(x) d X^{j}(x) \,,
\end{align*}
where $\Quad_{ij,k}$ and $\Quad_{ij}$ are smooth functions
on $\left.BM\right|_{\overline{W}}$.

Suppose that $n=1$.  Recall that $\DQ(\piQ\phi)$ is defined as
\begin{equation*}
	\DQ \piQ\phi  :=  d(\phi + f \n)    =  d\phi + \dQ f \wedge \n + f \wedge
	d\n  \,.
\end{equation*}
where $f\in\Wk{s-1}(M)$ is the unique function with $\n\wedge (d\phi +
f d\n)=0$. Because the map $h \mapsto h \n\wedge d\n$ is a smooth linear
isomorphism,  $\|f\|_{s-1} \prec \| f \n\wedge d\n\|_{s-1}$.
We can, therefore, estimate as follows:
\begin{equation*}
	\|\DQ \piQ\phi\|_{s-2}  \prec   \|d\phi\|_{s-2} +
	\|f\|_{s-1}   \prec   \|d\phi\|_{s-2} +
	\|d\phi \wedge \n\|_{s-1} \,.
\end{equation*}
 Thus, we need only estimate $\|d\phi\|_{s-2}$ and $\|d\phi \wedge \n\|_{s-1}$:
\begin{align*} 
\|d\phi\|_{W,s-2}  
=&\;
   \|d(\Quad_{ij,k}(x,X)X^{i}X^{j}\,dx^{k}
	+ \Quad_{ij}(x,X)X^{i}\,dX^{j})\|_{W,s-2}
\\
\prec&\;  \|X\|_{W,s-1}\|X\|_{W,s}
\\
\|(d\phi \wedge \n)\|_{W,s-1}   
\prec&\;
\|d (\Quad_{ij,k}(x,X)X^{i}X^{j}dx^{k}
	+ \Quad_{ij}(x,X)X^{i}dX^{j}) \wedge \n\|_{W,s-1}\\
\prec&\;  \|X\|_{W,s-1}\|X\|_{W,s}  \, .
\end{align*}
The proof of the second inequality follows by similar reasoning.

\medskip
Now suppose that $n>1$. Then
\begin{align*}
\|\dQ\piQ (\FX^* \n - \n - \Lie_{X})\|_{W,s-1}  
=&\; \| \piQ d \phi \|_{W,s-1}\\
=&\;
\| \dQ(\Quad_{ij,k}(x,X)X^{i}X^{j}\,dx^{k})
+ \dQ(\Quad_{ij}(x,X)X^{i}) \wedge \,\dQ X^{j})\|_{W,s-1}\\
\prec&\; \| X \|_{W,s} \, \|X \|_{W,s-1}
\end{align*}
The proof of the last  inequality in the statement of the lemma follows
by similar reasoning.
\end{proof}

\begin{proof}[Proof of Proposition~\ref{prop:Psi-estimates}]
Recall the definition of $\Psi$,
\begin{equation}
\label{eqn:Psi-1}
\Psi(X) - X  =  \Phi^{-1} \circ \pi \circ \Phi (X) - X
 =   \Phi^{-1} \circ \pi \circ \Phi (X) - \Phi^{-1}  \circ \Phi (X)
\end{equation}
where $\Psi$ is the map defined in Theorem~\ref{prop:localcoordn} and where
$\pi(g \oplus \a\oplus\omega) = g$.  Because $\Phi^{-1}$ is smooth (and
thus of class
$C^1$), Equation~\eqref{eqn:Psi-1} implies the inequality
\[
   \| \Psi(X) - X \|_{s}  \prec ||| \pi \circ \Phi(X) - \Phi(X) |||_{s}
	       = \| \piperp\Phi(X) \|_{s} \,,
\]
where $\piperp(g \oplus \a\oplus\omega) := \a \oplus\omega$ and
$||| g \oplus \a\oplus\omega |||_{s} :=
       \|g\|_{s+1} + \|\a \oplus  \omega\|_{s}$.
Consequently, to prove the estimate~\eqref{sharper-estimate}, we need only
estimate the terms in the expansion

\begin{multline}
\label{eqn:proj-phi}  \piperp\Phi(X) = 
\begin{cases}
(n+1) \GQ\deltaQ (\piQ \FX^* \n)  \oplus
	    (\GQ \DQs \DQ (\piQ\FX^* \n) \oplus  \HQ (\piQ \FX^* \n) ) \,,
 &\text{for $n=1$}\\
(n+1)\GQ \deltaQ (\piQ \FX^* \n)  \oplus
	     (n \GQ \deltaQ \dQ (\piQ\FX^* \n) \oplus  \HQ (\piQ \FX^* \n) )\,,
 &\text{for $n>1$}\\
\end{cases}
\end{multline}
given in Theorem~\ref{prop:localcoordn}.
By Lemma~\ref{lem:exp-q-form} we have the local formula
\[
\FX^{*} \n = \n + \Lie_{X} \n +
\Quad^1_{ij}(x,X)X^{i}X^{j} +\Quad^2_{ij}(x,X)X^{i}dX^{j} \,.
\]
Therefore, because  $X$ is contact and so $\piQ (\Lie_{X}\n) = 0$,
Lemma~\ref{lem:remainder-estimate} implies the estimate
\[
\| \piQ(\FX^* \n )\|_{s-1}
 \prec
    \| \piQ\n \|_{s-1} + \|\piQ (\Lie_{X}\n ) \|_{s-1}
      + \|X\|_s \, \|X\|_{s-1} = \|X\|_{s} \, \|X\|_{s-1} \,.
\]
Taking into account the orders of the various operators in
Equation~\eqref{eqn:proj-phi} and recalling
Lemma~\ref{lem:DQ-estimate},
we can estimate as follows:
\begin{align*}
\intertext{For $n\geq 1$,}
  \| \GQ\deltaQ \piQ (\FX^* \n) \|_{s} & \prec 
      \| \piQ (\FX^* \n) \|_{s-1} \prec \|X\|_{s}\,\|X\|_{s-1}
\\
  \|\HQ \piQ (\FX^* \n)\|_{s}& \prec  \| \piQ (\FX^*
	\n)\|_{s-1}   \prec \|X\|_{s}\,\|X\|_{s-1}\,;
\\
\intertext{for $n=1$,}
  \| \GQ \DQs \DQ \piQ(\FX^* \n) \|_{s} & \prec 
      \| \DQ \piQ(\FX^* \n) \|_{s-2}  \prec \|X\|_{s}\,\|X\|_{s-1}\,;
\intertext{and for $n>1$,}
  \| \GQ \deltaQ \dQ \piQ(\FX^* \n) \|_{s} & \prec 
      \| \dQ \piQ(\FX^* \n) \|_{s-1}  \prec \|X\|_{s}\,\|X\|_{s-1}\,.
\end{align*}
This concludes the proof of \eqref{sharper-estimate}.

We now show that
\begin{align*}
\|\left( \Psi(X_{2}) - X_{2} \right) - \left( \Psi(X_{1}) - X_{1} \right)\|_{s}
\prec  &\; \|X_{2} - X_{1}\|_{s-1}(\|X_{2}\|_{s} + \|X_{1}\|_{s})\\
 &\; + \|X_{2} - X_{1}\|_{s}(\|X_{2}\|_{s-1} + \|X_{1}\|_{s-1}) \,.
\end{align*}
Set $Y_i := \Psi(X_i) - X_i$, $i=1,2$. Differentiability of the map
 $\Phi^{-1}$ implies  the inequality
\[
  \|Y_2 - Y_1 \|_{s} =  \|\left( \Psi(X_{2}) - X_{2} \right)
	  - \left( \Psi(X_{1}) - X_{1} \right)\|_{s}
 \prec
    \|| \Phi( Y_2 - Y_1) \||_{s} \,;
\]
which, by adding and subtracting a term and applying 
the triangle inequality, implies the inequality
\begin{multline}
\label{eq:cauchy-est}
 \|Y_2 - Y_1 \|_s \prec
 ||| (\piperp\circ\Phi(X_2) - \piperp\circ\Phi(X_1) )|||_{s}\\
+    ||| \Phi( Y_2 - Y_1)  + (\piperp\circ\Phi(X_2) - \piperp\circ\Phi(X_1) )|||_{s}  \,.
\end{multline}
To conclude the proof, we need only estimate each term on the
right-hand side of \eqref{eq:cauchy-est}.

To estimate the first term, 
note that by Lemma~\ref{lem:exp-q-form} applied to the contact 
vector fields $X_i$,
we have  the local formula
\[
     \piQ F_{X_i}^*\n = \piQ(\n +  \Lie_{X_i}\n + \Quad_{\n}(X_i) )
       = \piQ (\Quad_{\n}(X_i) ) \,.
\]
This and Equation~\eqref{eqn:proj-phi}  imply the inequalities
\begin{align*}
 \|| \piperp\circ\Phi(X_2) - \piperp\circ\Phi(X_1) \||_{s}  
\prec
&\;
 \|\piQ( \Quad_{\n}(X_2)  - \Quad_{\n}(X_1) )\|_{s-1}\\
&\;     
+
 \|\DQ \piQ 
  \left( \Quad_{\n}(X_2) - \Quad_{\n}(X_1)
\right)\|_{s-2}\,
\text{ for $n=1$}\\
\intertext{and}
\prec
&\;
 \|\piQ( \Quad_{\n}(X_2)  - \Quad_{\n}(X_1) ) \|_{s-1}\\
&\;     
+
 \|\dQ \piQ 
  \left( \Quad_{\n}(X_2) - \Quad_{\n}(X_1)
\right)\|_{s-1}\,
\text{ for $n>1$}\,.
\end{align*}
Evoking  Lemma~\ref{lem:remainder-estimate} and  Lemma~\ref{lem:DQ-estimate}
then yields the estimate of the first term we need:
\begin{align*}
\label{eqn:cauchy-1}
\|| (\piperp\circ\Phi(X_2) - \piperp\circ\Phi(X_1) )\||_{s}   \prec &\;
       \|(X_{2}-X_{1})\|_{s-1}(\|X_{2}\|_{s}+\|X_{1}\|_{s})\\
&\; +  \|(X_{2}-X_{1})\|_{s}(\|X_{2}\|_{s-1}+\|X_{1}\|_{s-1})  \,.
\end{align*}

We estimate the second term in  \eqref{eq:cauchy-est} as follows.
We first claim that $\Phi(Y_2-Y_1)= \piperp\circ\Phi(Y_2-Y_1)$.
To see this, use the definition $\Psi(X) = \Phi^{-1}\circ\pi\circ\Phi(X)$ for any
vector field $X$ and the linearity of the map  $\pi\circ\Phi$ to write
$\pi\circ\Phi(Y_i)$ in the form
\[
\pi\circ\Phi(Y_i) = \pi\circ\Phi(\Psi(X_i)-X_i) = 
\pi\circ\Phi (\Phi^{-1}\circ\pi\circ\Phi(X_i)) - \pi\circ\Phi(X_i) =0\,.
\]
It follows that   $\pi\circ \Phi(Y_2-Y_1)=0$; hence,
$\piperp\circ\Phi(Y_2-Y_1) = \Phi(Y_2-Y_1)$.  
Next, since $\Psi(X) = \Phi^{-1} \circ\pi\circ\Phi(X)$, we have
\[
        \piperp\circ\Phi\circ\Psi(X_i) = \piperp\pi\circ\Phi(X_i)=0\,.
\] 
Combining these two identities shows the second term in
\eqref{eq:cauchy-est} can be written
\begin{multline*}
  \Phi( Y_2 - Y_1) 
   + (\piperp\circ\Phi(X_2) - \piperp\circ\Phi(X_1)  )\\
=    \piperp\left\{ \Phi( Y_2 - Y_1) 
   + ( \Phi(X_2) - \Phi(\Psi(X_2) )
      - \Phi(X_1) + \Phi(\Psi(X_1) )  \right\} .
\end{multline*}
We now show that
\begin{multline}
\label{eq:term2-b}
 \piperp\left\{ \Phi( Y_2 - Y_1) 
   + ( \Phi(X_2) - \Phi(\Psi(X_2) )
      - \Phi(X_1) + \Phi(\Psi(X_1) )
\right\} \\
 =
    L \left(\Quad_{\n}(Y_2-Y_1) -
	    (\Quad_{\n}( \Psi(X_2)) - \Quad_{\n}(\Psi(X_1)))
     + (\Quad_{\n}(X_2) - \Quad_{\n}(X_1))
       \right)
\end{multline}
where $L$ is the linear operator defined by
\[
   L(\phi):=
\begin{cases}
     \left\{(n+1)\, \GQ\deltaQ  \oplus (\GQ\DQs\DQ +\HQ ) \right\}\circ \piQ(\phi)  \,,
&\text{ for $n=1$}\\
     \left\{(n+1)\, \GQ \deltaQ   \oplus (n\,\GQ\deltaQ \dQ +\HQ ) \right\}\circ \piQ(\phi) \,,
&\text{ for $n>1$\,.}\\
\end{cases}
\]
To see this, apply Lemma~\ref{lem:exp-q-form} with $\phi=\n$ to the vector fields
$Y_2-Y_1$, $X_i$ and $\Psi(X_i)$ to obtain the three expansions 
\begin{align*}
F_{X_i}^*\n &= \n + \Lie_{X_i}\n + \Quad_{\n}(X_i)\\
F_{\Psi(X_i)}^*\n &= \n + \Lie_{\Psi(X_i)}\n 
   + \Quad_{\n}( \Psi(X_i))\\
F_{(Y_2 - Y_1)}^*\n &=  \Lie_{Y_2}\n -\Lie_{Y_1}\n
	+ \Quad_{\n}(Y_2 - Y_1) \,.
\end{align*}
Substituting these expressions into Equation~\eqref{eqn:proj-phi} 
and collecting terms reveals that the terms involving the Lie derivative of
$\n$ cancel to yield the identity \eqref{eq:term2-b}.
Thus, by the triangle inequality, to estimate the second term in
\eqref{eq:cauchy-est}, we need only estimate each term in the sum
\[
  ||| L ( \Quad_{\n}  (Y_2-Y_1) )  |||_{s}
+
 ||| L( \Quad_{\n} ( \Psi(X_2)) - \Quad_{\n}(\Psi(X_1))  )  |||_{s}
+
 ||| L( \Quad_{\n} (X_2) - \Quad_{\n}(X_1)  ) |||_{s}
\]
For each term, we use the estimates
\[
   |||L(\phi)|||_{s}\prec
\begin{cases}
   \|\piQ\phi\|_{s-1} + \| \DQ\piQ\phi\|_{s-2} + \|\HQ\piQ\phi\|_{s}  \,,
&\text{ for $n=1$}\\
    \|\piQ\phi\|_{s-1} + \| \dQ\piQ\phi\|_{s-1} + \|\HQ\piQ\phi\|_{s}  \,,
&\text{ for $n>1$\,.}\\
\end{cases}
\]
Noting the degree of the various linear operators in the definition of $L$
and employing  Lemmas~\ref{lem:remainder-estimate} and ~\ref{lem:DQ-estimate}
shows that
\begin{align*}
   ||| L (\Quad_{\n}  (Y_2-Y_1) ) |||_{s}
\prec &
(\| Y_2 - Y_1\|_{s} )\| Y_2 - Y_1 \|_{s-1} \,,\\
   ||| L (\Quad_{\n}(\Psi(X_2)) -  \Quad_{\n}(\Psi(X_1)) ) |||_{s}
\prec &
(\|\Psi(X_1)\|_s + \|\Psi(X_2)\|_{s} ) \cdot \|\Psi(X_2) - \Psi(X_1) \|_{s-1}\\
&\; + (\|\Psi(X_1)\|_{s-1} + \|\Psi(X_2)\|_{s-1} ) \cdot \|\Psi(X_2) - \Psi(X_1) \|_{s}\,,
\intertext{and}
   ||| L (\Quad_{\n}(X_2) -  \Quad_{\n}(X_1) ) |||_{s}
\prec &
(\| X_1 \|_s + \| X_2 \|_{s} ) \cdot \|X_2 - X_1 \|_{s-1} \\
&\; + (\| X_1 \|_{s-1} + \| X_2 \|_{s-1} ) \cdot \|X_2 - X_1 \|_{s}\\
\end{align*}
To conclude the estimate of \eqref{eq:cauchy-est}, note that since $\Psi$
is $C^1$, we can replace $\Psi(X_{i})$ by $X_{i}$ and $Y_2-Y_1$ by
$X_2-X_1$ everywhere in the first two of the previous three inequalities.
\end{proof}



\begin{thebibliography}{9999}


%

\bibitem[Arn]{Arn:1978}
	 V. I. Arnold,
	 {\em Mathematical Methods of Classical Mechanics\/},
	 Springer-Verlag
	 (1978),
	 New York, Heidelberg, Berlin.

%
\bibitem[B]{B}
	  J.~Bland,
	  {\em Contact geometry and {CR}-structures on {$S^3$}},
	   Acta Math.
	  {\bf 172} (1994),
	  \mbox{1--49}.
 \bibitem[BD2]{BD2}
	 J. Bland and T. Duchamp,
	 {\em Anisotropic Estimates for Sub-elliptic  Operators \/},
	  Sci. in China, Ser. A: Mathematics {\bf 51} 
	 (Science Press/ Springer-Verlag) (2008), \mbox{509--522}.  
%
 \bibitem[BD3]{BD3}
	 J. Bland and T. Duchamp,
	 {\em The space of Cauchy-Riemann structures on 3-D compact contact manifolds \/},
	(preprint).  
%
\bibitem[Bi]{Bi}
	  O. Biquard,
	   {\em Metriques autoduales sur la boule\/},
	   Invent. math.
	  {\bf 148}~(2002),
	  \mbox{545-607}.
%
 \bibitem[CL]{CL}
	  J.-H. Cheng and J.M. Lee,
	  {\em A local slice theorem for 3-dimensional CR
	  structures\/},
	  Amer. J. Math. {\bf 117} (1995),
	  \mbox{1249--1298}.
%
\bibitem[E]{Ebin:1970}
	D. Ebin, {\em The manifold of Riemannian metrics},
	in ``Global Analysis'', eds. S.S. Chern and S. Smale, Proc.
	Symp. in Pure Math., Amer. Math. Soc. Vol XV (1970),
	11--40.
%
\bibitem[EM]{Ebin-Marsden:1970}
	 D.~Ebin and J.~Marsden,
	{\em Groups of diffeomorphisms and the motion of an incompressible fluid},
	Ann. Math. \textbf{92} (1070)
	\mbox{102--163}.
%
 \bibitem[FS]{FS}
	    G. B. Folland and E. M. Stein,
	 {\em Estimates for the $\db_{b}$ complex and analysis on
	 the Heisenberg group\/},
	 Comm. on Pure and Applied Math.
	 {\bf 27}~(1974),
	 \mbox{429--522}.
%
\bibitem[G]{gray:1959a}
	J.~Gray, {\em Some global properties of contact structures}, Ann. of Math.
	  {\bf 69} (1959), 421--450.
%
\bibitem[Ham]{Ham}
	 R. Hamilton,
	 {\em The inverse function theorem of Nash and Moser},
	 Bull. AMS
	 \textbf{7} (1982)
	 \mbox{65--222}.
%
%
\bibitem[O1]{Omori:1970}
H.~Omori, {\em On the group of diffeomorphisms on a compact manifold},
in ``Global Analysis'', eds. S.S. Chern and S. Smale, Proc.
Symp. in Pure Math., Amer. Math. Soc. Vol XV (1970), 167--183.
%
\bibitem[O2]{Omori:1974}
 H.~Omori,
 {\em Infinite dimensional lie transformations groups},
 Lecture Notes in Mathematics, vol. 427,
 Springer-Verlag, Berlin, Heidelberg, New York, 1974.
 %
 \bibitem[Pal]{Palais:1968}
	 R.  Palais,
	 {\em Foundations of Global Non-Linear Analysis\/},
	 W.A. Benjamin, Inc.
	 (1968),
	 New York.
 %
\bibitem[R]{Rumin:1994}
      M.~Rumin,
{\em Formes differentielles sur les varietes de contact},
J. Diff. Geom. \textbf{39} (1994), \mbox{281--330}.

\end{thebibliography}
\end{document}